\documentclass[11pt]{amsart}
\usepackage{amssymb,amsmath,txfonts,mathrsfs,mathtools,titletoc, tikz, stmaryrd}
\usepackage[alphabetic]{amsrefs}
\usepackage{amsthm}
\newtheorem{theorem}{Theorem}[section]
\newtheorem{prop}[theorem]{Proposition}
\newtheorem{lemma}[theorem]{Lemma}
\newtheorem{remark}[theorem]{Remark}

\newtheorem{definition}[theorem]{Definition}
\newtheorem{cor}[theorem]{Corollary}

\newtheorem{conj}[theorem]{Conjecture}

\usepackage{float}
\usepackage[colorlinks, linkcolor=blue,anchorcolor=Periwinkle,
    citecolor=red,urlcolor=Emerald]{hyperref}
\usepackage{xcolor}
\usepackage{enumitem}
\usepackage{geometry,array} 
\usepackage{graphicx}
\graphicspath{{./graphic/}}
\usepackage{subfigure}
\usepackage{bookmark}
\usepackage{tikz}\usetikzlibrary{matrix}
\usepackage{url}
\usepackage{braids}

\usetikzlibrary{decorations.markings}
\tikzset{->-/.style={decoration={  markings,  mark=at position #1 with
    {\arrow{>}}},postaction={decorate}}}
\tikzset{-<-/.style={decoration={  markings,  mark=at position #1 with
    {\arrow{<}}},postaction={decorate}}}
\usepackage[all]{xy}
\usetikzlibrary{arrows}
\usetikzlibrary{decorations.pathreplacing,decorations.pathmorphing,shadings,fadings,calc}
\usepackage{extarrows}
\usepackage{tikz-3dplot}

\def\pf{{\it Proof:}~}

\usetikzlibrary{matrix}
\usepackage[all]{xy}
\usetikzlibrary{arrows,shapes}
\usetikzlibrary{calc}
\tikzset{->-/.style={decoration={  markings,  mark=at position #1 with
    {\arrow{>}}},postaction={decorate}}}
\tikzset{-<-/.style={decoration={  markings,  mark=at position #1 with
    {\arrow{<}}},postaction={decorate}}}
\usepackage[]{hyperref}
\numberwithin{equation}{section}

\usetikzlibrary{decorations.markings}
\tikzset{->-/.style={decoration={  markings,  mark=at position #1 with
    {\arrow{>}}},postaction={decorate}}}
\tikzset{-<-/.style={decoration={  markings,  mark=at position #1 with
    {\arrow{<}}},postaction={decorate}}}

\usepackage{lscape}
\usepackage{CJK}

\numberwithin{equation}{section}

\begin{document}
\title{Multiple eigenvalues and the width}
\author{Qixuan Hu}
\address{Qixuan Hu\\ Department of Mathematics and Computer Sciences\\ Shantou University, Guangdong\\P. R. China}
\email{qxhu@stu.edu.cn}
\date{\today}

\begin{abstract}
We obtain the simplicity of the first Neumann eigenvalue of convex `thin' domain with boundary in $\mathbb{R}^n$ and compact `thin' manifolds with non-negative Ricci curvature. For convex 'thin' domain in $\mathbb{R}^2$, we get the simplicity of the first $k$ Neumann eigenvalues. The number $k$ depends on the ratio of the corresponding width over the diameter of the domain.  For convex 'thin' domain in $\mathbb{R}^n$, we obtain the eigenvalue comparison with collapsing segment.
\\[3mm]
Mathematics Subject Classification: 35K15, 53C44
\end{abstract}
\thanks{Q. Hu was partially supported by SRIG NTF25026T.}

\maketitle

\titlecontents{section}[0em]{}{\hspace{.5em}}{}{\titlerule*[1pc]{.}\contentspage}
\titlecontents{subsection}[1.5em]{}{\hspace{.5em}}{}{\titlerule*[1pc]{.}\contentspage}
\tableofcontents

\section{Introduction}

Based on the simplicity of the first several Neumann eigenvalues of thin rectangle (which is a convex domain with boundary) and thin tube (which is a convex surface without boundary),  we study the simplicity of Neumann eigenvalues of convex `thin' domain with boundary in $\mathbb{R}^n$ and compact `thin' manifolds with non-negative Ricci curvature .  

In $1976$,  Cheng \cite{Cheng} showed that for a compact Riemannian surface $M^2$ of genus $g$, the multiplicity of $\mu_i(M^2)$ is $\leq \frac{(2g+ i+ 1)(2g+ i+ 2)}{2}$.  Later, Besson \cite{Besson} improved this quadratic bound to a linear bound $4g+2i+1$. In \cite{Na}, Nadirashvili improved Besson's result even further. Also see \cite{HHN} and \cite{HMN} for further discussion.  In \cite{Se} S\'{e}vennec proved that the first eigenvalue has multiplicity at most $5-\chi(M)$ if $\chi(M)<0$.

Also Nadirashvili \cite{Na} proved that the multiplicity of $\mu_1$ is at most 2.  It is well-known that for a round ball or a square in $\mathbb{R}^2$, their first Neumann eigenvalue has multiplicity 2. 

A natural question to ask is under what kind of restriction on the domain can we get the first eigenvalue is simple.  

\begin{remark}
	The discussion for bounds on multiplicity of eigenvalues by constrain on the topology of manifold focused on 2 dimensional case, we cannot expect similar result for higher dimension due to the following theorem of Colin de Verdi\`{e}re. 
\end{remark}

\begin{theorem}\cite{CdV}\label{high dimension no bound}
	Let $M$ be a compact manifold of dimension $\geq3$. Then for any finite sequence $0=a_0<a_1\leq\cdots\leq a_n$, there exists a $C^\infty$ metric $g$ on $M$, such that this sequence is the first $n+1$ eigenvalue of the Laplace operator on $(M,g)$.
\end{theorem}

\begin{definition}\label{def width of domain}
	{For $\Omega\subseteq \mathbb{R}^n$ be a bounded convex domain, let $\mathcal{D}(\Omega)$ be the diameter of $\Omega$ in $\mathbb{R}^n$, we define $\mathcal{W}(\Omega)$ the width of $\Omega\subseteq \mathbb{R}^n$ as follows:
		\begin{align}
			&\mathscr{P}_{v}(\Omega)\vcentcolon= \{x\in \mathbb{R}^n: x\cdot v= 0, x= y+ s\cdot v \ \text{for some}\ y\in \Omega, s\in \mathbb{R}\}, \quad  \forall v\in \mathbb{S}^{n- 1}\subseteq \mathbb{R}^n; \nonumber \\
			&\mathcal{W}(\Omega)= \min_{v\in \mathbb{S}^{n- 1}}  \mathcal{D}(\mathscr{P}_{v}(\Omega)). \nonumber
		\end{align}
	}
\end{definition}

Ba\~{n}uelos and Burdzy \cite{BB} proved the following result,
\begin{theorem}\cite{BB}
	Let $\Omega\subset\mathbb{R}^2$ be a convex domain. If the ratio of the width and diameter of $\Omega$ is small, $$\frac{\mathcal{W}(\Omega)}{\mathcal{D}(\Omega)}\leq \frac{\pi}{4j_{0,1}},$$
	where $j_{0,1}$ is the first zero of zero order Bessel function $J_0(x)$. Then the first Neumann eigenvalue of $\Omega$, $\mu_1$ is simple.
\end{theorem}

 Their proof use probabilistic method.  They also made the following conjecture:
\begin{conj}\label{conj square simple critical}
{If $\mu_1(\Omega)$ has multiplicity $=2$ for some convex domain $\Omega\subseteq \mathbb{R}^2$,  then 
\begin{align}
\frac{\mathcal{W}(\Omega)}{\mathcal{D}(\Omega)}\geq \frac{1}{\sqrt{2}}.  \nonumber 
\end{align}
}
\end{conj}

In this paper, we give an analytic proof which extends Ba\~{n}uelos and Burdzy's result to higher dimension and higher order of eigenvalue.

In the first part of this paper, we prove the following statement:
\begin{theorem}\label{main statement 1}
{For any positive integer $k$. There exists a constant $\epsilon_k=\frac{\pi}{4(j_{0,1}+\frac{k-1}{2}\pi)}$, such that for any convex domain $\Omega\subset\mathbb{R}^2$ with $C^1$ smooth boundary,  if $\frac{\mathcal{W}(\Omega)}{\mathcal{D}(\Omega)}<\epsilon_k$, then the first $k$ Neumann eigenvalue of $\Omega$ is simple.
}
\end{theorem}

\begin{theorem}\label{main statement 2}
	{For any positive integer $n$. There exists a constant $\delta_n=\frac{\sqrt{C_n}}{\sqrt{8}j_{\frac{n-2}{2},1}}$, where $C_n=-\Gamma(\frac{n}{2})(\frac{2}{x_m})^{\frac{n-2}{2}}J_{\frac{n-2}{2}}(x_m)$, and $x_m$ is the first minimum point of $t^{-\frac{n-2}{2}}J_{\frac{n-2}{2}}(t)$. such that for any convex domain $\Omega\subset\mathbb{R}^n$ with $C^1$ smooth boundary,  if $\frac{\mathcal{W}(\Omega)}{\mathcal{D}(\Omega)}<\delta_n$, then the first Neumann eigenvalue of $\Omega$ is simple.
	}
\end{theorem}

\begin{theorem}\label{main statement 3}
	For a compact complete $n$-dim Riemannian manifold $M$ with $Rc\geq0$. Let $A,B\in M$, such that $d(A,B)=\mathcal{D}(M)$, let $l$ be a segment in $M$ connecting $A,B$. Suppose that $M$ is in the $\epsilon$ neighborhood of $l$. If the first eigenvalue $\mu$ of $M$ is not simple, then we have
	$$
	\mu\geq\frac{C(n)}{\epsilon^2},
	$$
	where $C(n)=\frac{C_n}{18}$, for $C_n$ defined in Theorem \ref{main statement 2}.
\end{theorem}

In the second part of this paper, we compare the eigenvalue of domain $\Omega$ with the eigenvalue of the collapsing segment. 

Let $\Omega\subset\mathbb{R}^n$ be a convex domain, and $\mathcal{D}(\Omega)=1,\mathcal{W}(\Omega)=\epsilon$.

Assume $\Omega^t=\{y\in\mathbb{R}^{n-1}:(t,y)\in\Omega\}$ and $\Omega=\{(t,y):t\in[0,1], y\in\Omega^t\}$.
We define
\begin{align}
	H(x)\vcentcolon= \frac{\mathcal{H}^{n-1}(\Omega^x)}{\int_0^1 \mathcal{H}^{n-1}(\Omega^t)dt},  \nonumber 
\end{align}
where $\mathcal{H}^{n- 1}$ is $(n-1)$-dim Hausdorff measure.

Let $d\mu(x)=H(x)dx$, and let $N$ be the segment $[0,1]$ with measure $d\mu$. 

For convex domain in $\mathbb{R}^2$, Wang and Xu\cite{WX} proved the following result,
\begin{theorem}\cite{WX}
	There is a universal constant $C_1,C_2>0$, such that for any $\Omega\subset\mathbb{R}^2$ be a convex domain with diameter $\mathcal{D}(\Omega)=1$, and width $\mathcal{W}(\Omega)=\epsilon\in(0,C_1)$, we have 
	$$
	\mu_1(\Omega)\geq\pi^2+C_2\epsilon^2.
	$$
\end{theorem}

In addition, Wang and Xu \cite{WX}[Prop 6.3] show that there exists a universal constant $C>0$, such that for convex domain $\Omega\subset\mathbb{R}^n$, we have 
$$
\mu_1(N)-C\epsilon\leq\mu_1(\Omega)\leq\mu_1(N).
$$
This result is due to \cite{Jer} for $n=2$, and \cite{CJK} for $n\geq2$.
Here, we give a different proof which also applies to higher order eigenvalue.

\begin{theorem}\label{main statement 4}
Let $\mu_k(N),\mu_k(\Omega)$ be the $k$-th Neumann eigenvalue of $N,\Omega$. We have the following comparison,
\begin{equation}
(1-2\epsilon(1+\mu_k(N)))\cdot\mu_k(N)\leq\mu_k(\Omega)\leq\mu_k(N). \nonumber
\end{equation}
\end{theorem}

Using similar argument as in \cite{WX},  we can improve the inequality in Theorem \ref{main statement 4} for the case  $\Omega$ is a convex domain in $\mathbb{R}^2$, and $k=1$, we have the following comparison,

\begin{theorem}
	There is a universal constant $C>0$, for any $\Omega\subset\mathbb{R}^2$ be a convex domain, we have
\begin{equation*}
	\mu_1(N)-C\epsilon^2\leq\mu_1(\Omega)\leq\mu_1(N).
\end{equation*}
\end{theorem}

\part{The Multiplicity of Neumann Eigenvalues}

\section{The multiple eigenvalue of convex planar domain}

In this section, assume $\Omega$ is a convex domain in $\mathbb{R}^2$ with $C^1$ smooth boundary,  and width $\mathcal{W}(\Omega)=\mathcal{D}(\mathscr{P}_{v_0}(\Omega))= \epsilon$, where $v_0=(1,0)$. 

\begin{lemma}\label{lem upper and lower part of boundary}
	{Choose $z_L=(z^1_L,z^2_L), z_R=(z_R^1,z_R^2)\in\overline{\Omega}$ such that 
		\begin{align}
			z_L^1= \min\{x\in \mathbb{R}: (x, y)\in \Omega\}, \quad \quad \quad z_R^1= \max\{x\in \mathbb{R}: (x, y)\in \Omega\}. \nonumber 
		\end{align}
		
		then $z_L$ and $z_R$ divides $\partial\Omega$ into two parts $\Gamma_1$ and $\Gamma_2$, such that

		\begin{align}
			\vec{n}\big|_{\Gamma_1}\cdot \vec{e}_2\geq 0, \quad \quad \vec{n}\big|_{\Gamma_2}\cdot \vec{e}_2\leq 0, \nonumber 
		\end{align}
		where $\vec{e}_1=(1,0)$ and $\vec{e}_2= (0, 1)\in \mathbb{R}^2$, $\vec{n}\big|_{\Gamma_i}$ is the outward normal vector of $\Omega$ on $\Gamma_i$ with $i= 1, 2$.  
	}
\end{lemma}

\pf
Since $\Omega$ is convex, we get the segment $\overline{z_Lz_R}$ divides $\Omega$ into two parts $D_1$, $D_2$, where 
$$D_1=\{z\in\Omega, \mathrm{the\, ratio\, of\, line}\, z_Lz\, \mathrm{is\, larger\, than\, the\, ratio\, of\, line}\, z_Lz_R \}
$$ 
and 
$$D_2=\{z\in\Omega, \mathrm{the\, ratio\, of\, line}\, z_Lz\, \mathrm{is\, smaller\, than\, the\, ratio\, of\, line}\, z_Lz_R \}
$$
Let $\Gamma_1=\partial \Omega\cap \overline{D_1}$, and
$\Gamma_2=\partial\Omega\cap\overline{D_2}$,

For any $p=(x,y)\in\Gamma_1$, we have $z_L^1\leq x\leq z_R^1$, if $x=z_L^1$, from convexity of $\Omega$ and $z_L^1= \min\{x\in \mathbb{R}: (x, y)\in \Omega\}$, we get segment $\overline{z_Lp}\subset\Gamma_1$ and $\vec{n}\Big|_{p}=(-1,0)$, so $\vec{n}\cdot \vec{e}_2=0$ at $p$. Similarly, if $x=z_R^1$, we get $\vec{n}\Big|_p=(1,0)$ and $\vec{n}\cdot \vec{e}_2=0$ at $p$.

If $x\in(z_L^1,z_R^1)$, let $\vec{t}=e^{i\theta}$ be the unit tangent vector of $\partial\Omega$ at $p$, where $\theta\in(-\pi,\pi)$ such that $\vec{n}=\vec{t}\cdot e^{\frac{i\pi}{2}}$, by convexity of $\Omega$, we have $\tan\theta$ is a decreasing function of $x$, and $\lim_{x\to z_L^1}\theta(x)=\frac{\pi}{2}$, $\lim_{x\to z_R^1}\theta(x)=-\frac{\pi}{2}$, so for $x\in(z_L^1,z_R^1)$, $\theta(x)\in(-\frac{\pi}{2},\frac{\pi}{2})$, and $\vec{n}\cdot \vec{e}_2=\vec{t}\cdot \vec{e}_1=\cos\theta>0$.

Similar argument on $\Gamma_2$ shows $\vec{n}\cdot \vec{e}_2\leq0$ on $\Gamma_2$.
\qed

\begin{theorem}\label{thm non-simple eigenvalue lower bound in width}
{Let $\Omega\subset\mathbb{R}^2$ be a convex domain with $C^1$ smooth boundary, and the width $\mathcal{W}(\Omega)=\epsilon$. Assume $\mu$ is a multiple Neumann eigenvalue of $\Omega$, then $\displaystyle \mu\geq \frac{\pi^2}{4\epsilon^2}$. 
}
\end{theorem}

\pf
\textbf{Step (1)}. By rotation, we may assume $\mathcal{W}(\Omega)=\mathcal{D}(\mathscr{P}_{v_0}(\Omega))= \epsilon$, where $v_0=(1,0)$.

Let $u_1$, $u_2$ be two linear independent eigenfunction of $\Omega$ with respect to eigenvalue $\mu$.

If $u_1(z_L)=0$, we let $u=u_1$, otherwise, we let $u=u_1(z_L)u_2-u_2(z_L)u_1$. We have $u$ is also eigenfunction of $\Omega$ with respect to $\mu$, and $u(z_L)=0$. Let $\Sigma= u^{-1}(0)$. 

By Courant-Cheng's nodal domain theorem, we have $\Omega\setminus\Sigma=\cup_{i=1}^s\Omega^i$ where $s$ is finite and $\Omega^1,\cdots,\Omega^s$ be all the nodal domains of $u$. 

Note $\Sigma$ is one-dimensional manifold except on a zero-dimensional subset of $\Omega$, so $\Omega\setminus\Sigma$ is dense in $\Omega$. Therefore we get 
\begin{align}
\bar{\Omega}=\bigcup_{i=1}^s\bar{\Omega^i}, \quad \quad \partial\Omega\subset \bigcup_{i=1}^s\partial{\Omega^i}. \nonumber 
\end{align}

\textbf{Step (2)}. For $\{z_k\}_{k= 1}^\infty\subseteq \Gamma_1$ with $\displaystyle \lim_{k\rightarrow\infty} z_k= z_L$, we claim that there exists $N>0$ and unique $i\in\{1,\cdots,s\}$, such that $z_ k\in\partial\Omega^i$ for $k>N$. 

Otherwise, we can find $z_{k_1},\cdots,z_{k_4}$ which is arranged clockwise on $\Gamma_1\subset\partial\Omega$, and there exists $\alpha,\beta\in\{1,\cdots,s\}$, such that $\alpha\neq\beta$ and $z_{k_1}$, $z_{k_3}\in \partial\Omega^\alpha$ and $z_{k_2}$, $z_{k_4}\in \partial\Omega^\beta$. 

Choose simple curve $\gamma_1$ in $\Omega^\alpha$ connecting $z_{k_1},\,z_{k_3}$, and simple curve $\gamma_2$ in $\Omega^\beta$ connecting $z_{k_2},\,z_{k_4}$, by the arrangement of 4 points $z_{k_1},\cdots,z_{k_4}$ on $\partial\Omega$, we get $\gamma_1\cap\gamma_2\neq\emptyset$, it is the contradiction. 

Without loss of generality, there exists $\delta>0$, such that 
\begin{align}
\Big[B_{z_L}(\delta)\cap\Gamma_1\Big]\subset\partial\Omega^1, \quad \quad \Big[B_{z_L}(\delta)\cap\Gamma_1\Big]\cap\partial\Omega^k= \emptyset, \quad \quad \quad \forall k\neq 1. \label{ball and gamma i intersection}
\end{align}
We have similar conclusion for $\Gamma_2$.

Now we assume that for some $j\in\{1,\cdots,s\}$, 
\begin{align}
[B_{z_L}(\delta)\cap\Gamma_1]\subset\partial\Omega^1, \quad \quad [B_{z_L}(\delta)\cap\Gamma_2]\subset\partial\Omega^j. \nonumber 
\end{align}

Without loss of generality, there are only two cases as follows left:
\begin{align}
[B_{z_L}(\delta)\cap\Gamma_2]\subset\partial\Omega^2 \quad \text{or}\quad [B_{z_L}(\delta)\cap\Gamma_2]\subset\partial\Omega^1.  \nonumber 
\end{align}

\textbf{Step (3)}. If $[B_{z_L}(\delta)\cap\Gamma_2]\subset\partial\Omega^2$, we claim 
\begin{align}
(\partial\Omega^1\cap\partial \Omega)\subset \Gamma_1 \quad \text{or}\quad (\partial \Omega^2\cap\partial\Omega)\subset\Gamma_2.  \nonumber 
\end{align}

If $(\partial\Omega^1\cap\partial \Omega)\subset \Gamma_1$ does not hold, then we can find $P\in B_{z_L}(\delta)\cap\Gamma_1$ and $Q\in\partial \Omega^1\cap \Gamma_2$. 

Assume $L_1$ and $L_2$ are two components of $\partial\Omega\setminus\{P,Q\}$ with $z_L\in L_1$.  Let $\gamma_1$ be a simple curve in $\Omega^1$ connecting $P,Q$. Since $\Omega^2\cap\gamma_1= \emptyset$, the set $\Omega^2$ is contained in one component of $\Omega\setminus\gamma_1$.

Note $z_L\in L_1$ and $z_L\in\partial\Omega^2$, the set $\Omega^2$ is contained in the component of $\Omega\setminus\gamma_1$, which has $L_1$ and $\gamma_1$ as boundary. So $[\partial\Omega^2\cap\partial\Omega]\subset L_1$. 

By $[L_1\cap\Gamma_1]\subset [B_{z_L}(\delta)\cap\Gamma_1]$ and (\ref{ball and gamma i intersection}), we get 
\begin{align}
[\partial\Omega^2\cap\partial\Omega\cap \Gamma_1]\subseteq [\partial\Omega^2\cap L_1\cap  \Gamma_1]\subseteq [\partial\Omega^2\cap B_{z_L}(\delta)\cap\Gamma_1]= \emptyset. \nonumber 
\end{align}
Hence $\displaystyle [\partial\Omega^2\cap\partial\Omega]\subseteq \Gamma_2$, applying Theorem \ref{thm DN-heat kernel comparison} for $\Omega^2$, we get 
\begin{align}
\mu\geq \lambda_1(\Omega^2)\geq \frac{\pi^2}{4\epsilon^2}. \nonumber
\end{align}

\textbf{Step (4)}. Finally we deal with the case $[B_{z_L}(\delta)\cap\Gamma_2]\subset\partial\Omega^1$. We assume $u>0$ in $\Omega^1$. 

By Lemma \ref{lem positive negative} and there are only finite nodal domains, we can choose a sequence $z_k\to z_L$ such that $u(z_k)<0$, and $z_k\in \Omega^3$ for all $k$. And we get $z_L\in\partial\Omega^3$. 

We choose $P\in B_{z_L}(\delta)\cap\Gamma_1$, and $Q\in B_{z_L}(\delta)\cap\Gamma_2$. Let $\gamma_1$ be a simple curve in $\Omega^1$ connecting $P$ and $Q$. Let $L_1$ and $L_2$ be two components of $\partial\Omega\setminus\{P,Q\}$, and $L_1$ is the component which contain $z_L$. 

Since $z_L\in\partial\Omega^3$, we get $\Omega^3$ is contained in the domain enclosed by $\gamma_1\cup L_1$. So $\partial\Omega^3\cap\partial\Omega\subset L_1$. By (\ref{ball and gamma i intersection}), we have
\begin{align}
\big([\partial\Omega^3\cap\partial\Omega]\backslash \{z_L\}\big)\subseteq [\partial\Omega^3\cap B_{z_L}(\delta)\cap\Gamma_1]\cup [\partial\Omega^3\cap B_{z_L}(\delta)\cap\Gamma_2]= \emptyset .\nonumber 
\end{align}

Therefore we have $\displaystyle [\partial\Omega^3\cap\partial\Omega]=\{z_L\}$. Now applying Theorem \ref{thm DN-heat kernel comparison} for $\Omega^3$(the Dirichlet boundary is empty in this case), we get $\mu=\lambda_1(\Omega^3)\geq\frac{\pi^2}{4\epsilon^2}$.
\qed

We recall the following result from \cite{Kr}.
\begin{theorem}\label{upper bound of eigen of N}
	{Assume $\Omega\subset\mathbb{R}^n$ is a convex domain with diameter $\mathcal{D}(\Omega)=1$,let $k\in \mathbb{Z}^+$, and $\mu_k(\Omega)$ be the $k$-th Neumann eigenvalue of $\Omega$, then we have
		\begin{enumerate}
			\item[(a)]. if $n= 2$, then $\mu_k(\Omega)\leq 4(j_{0, 1}+ \frac{k- 1}{2}\pi)^2$;
			\item[(b)]. if $n> 2$, and $k$ is odd, then $\mu_k(\Omega)\leq 4j_{\frac{n-2}{2}, \frac{k+ 1}{2}}^2$;
			\item[(c)]. if $n> 2$, and $k$ is even, then $\mu_k(\Omega)\leq (j_{\frac{n-2}{2}, \frac{k}{2}}+ j_{\frac{n-2}{2}, \frac{k+ 2}{2}})^2$;
		\end{enumerate}
	}
\end{theorem}

From Theorem \ref{upper bound of eigen of N} and Theorem \ref{thm non-simple eigenvalue lower bound in width}, we get the following result,

\begin{theorem}\label{2dim eigen simple}
	Let $\Omega\subset\mathbb{R}^2$ be a convex domain with $C^1$ boundary, $\Omega$ has diameter $\mathcal{D}(\Omega)$ and width $\mathcal{W}(\Omega)$, if $\frac{\mathcal{W}(\Omega)}{\mathcal{D}(\Omega)}<\frac{\pi}{4(j_{0,1}+\frac{k-1}{2}\pi)}$, then the first $k$ Neumann eigenvalue of $\Omega$ is simple.
\end{theorem}

\pf
By scaling, we assume $\mathcal{D}(\Omega)=1$, and $\mathcal{W}(\Omega)=\epsilon$. Suppose $\mu$ is a multiple eigenvalue of $\Omega$, by Theorem \ref{thm non-simple eigenvalue lower bound in width}, we get $\mu\geq\frac{\pi^2}{4\epsilon^2}$. If $\epsilon<\frac{\pi}{4(j_{0,1}+\frac{k-1}{2}\pi)}$, we get $\mu>4(j_{0,1}+\frac{k-1}{2}\pi)^2$. By Theorem \ref{upper bound of eigen of N}, we get $\mu$ cannot be first $k$ Neumann eigenvalue of $\Omega$.
\qed

\section{The first Neumann eigenvalue of convex domain}

In this section, we assume $\Omega\subset\mathbb{R}^n$ is a convex domain with $C^1$ smooth boundary, and  $\mathcal{W}(\Omega)= \epsilon> 0$. Furthermore, we assume $\mathcal{D}(\mathscr{P}_{v_0}(\Omega))=\mathcal{W}(\Omega)=\epsilon$, where $v_0=(1,0,\cdots,0)$. From definition of $\mathscr{P}_v(\Omega)$, we get $\Omega\subset\mathbb{R}\times\mathscr{P}_{v_0}(\Omega)$.

\begin{lemma}\label{us0 vertical control}
	Let $u$ be  the first Neumann eigenfunction of $\Omega\subset\mathbb{R}^n$ with respect to eigenvalue $\mu$, such that $\sup u=1,\,\inf u=-k$, and $0<k\leq1$. Suppose $u(s_0,\mathbf{y}_0)=-k,\,u(s_1,\mathbf{y}_1)=1$, for some $(s_0,\mathbf{y}_0),(s_1,\mathbf{y}_1)\in\overline{\Omega}$, Then for any $(s_0,\mathbf{y})\in\overline{\Omega}$, we have $u(s_0,\mathbf{y})\leq-k+2\mu\epsilon^2$, and for any $(s_1,\mathbf{y})\in\overline\Omega$, we have $u(s_1,\mathbf{y})\geq1-2\mu\epsilon^2$.
\end{lemma}

\pf
Let $v=\frac{u-\frac{1-k}{2}}{\frac{1+k}{2}}$ and $\theta=\sin^{-1}v\in[-\frac{\pi}{2},\frac{\pi}{2}]$. From gradient estimate of eigenfunction \cite{ZY}, we have
\begin{align}
	|D\theta|\leq\sqrt{\mu}\sqrt{\frac{2}{1+k}}\leq\sqrt{2\mu}.
\end{align}
So, we get
\begin{align}\label{Du est}
	|Du|=\frac{1+k}{2}|Dv|=\frac{1+k}{2}\sqrt{1-v^2}|D\theta|\leq\sqrt{(u+k)(1-k)}\sqrt{2\mu}.
\end{align}

For any $(s_0,\mathbf{y})\in\Omega$, since $\Omega\subset\mathbb{R}\times\mathscr{P}_{v_0}(\Omega)$, and $\mathcal{D}(\mathscr{P}_{v_0}(\Omega))=\epsilon$, we get $|(s_0,\mathbf{y})-(s_0,\mathbf{y}_0)|\leq\epsilon$, so we have
\begin{align}\label{us0 est}
	u(s_0,\mathbf{y})\leq-k+\sup|Du|\cdot \epsilon\leq-k+\epsilon\sqrt{2\mu}.
\end{align}
Plug in (\ref{us0 est}) in (\ref{Du est}), we get
\begin{align}\label{Dus0 est}
	|Du|(s_0,y)\leq \sqrt{2\mu}\sqrt{u(s_0,y)+k}\leq\sqrt{2\mu}(\epsilon\sqrt{2\mu})^\frac{1}{2}.
\end{align}
Plug in (\ref{Dus0 est}) in (\ref{us0 est}) we get
\begin{align}\label{us0 rst 2}
	u(s_0,\mathbf{y})\leq-k+\epsilon\cdot\sqrt{2\mu}(\epsilon\sqrt{2\mu})^\frac{1}{2}.
\end{align}
Repeat this process $m$ times, by induction we get
\begin{align}
	|Du|(s_0,\mathbf{y})\leq\sqrt{2\mu}(\epsilon\sqrt{2\mu})^{2^{-1}+2^{-2}+\cdots+2^{-m}}.
\end{align}
Let $m\to\infty$, we get
\begin{align}
	|Du|(s_0,\mathbf{y})\leq 2\mu\epsilon,
\end{align}
and
\begin{align}
	u(s_0,\mathbf{y})\leq -k+2\mu\epsilon^2.
\end{align}

Similar arguments holds for $u(s_1,\mathbf{y})$ and $|Du|(s_1,\mathbf{y})$.
\qed

From \cite{Kr2}, we get the following estimate on the ratio of $\sup u$ and $-\inf u$. 

\begin{lemma}\label{max min value mfd}\cite{Kr2}
	Let $M$ be a compact $n$-dim Riemannian manifold with non-negative Ricci curvature, or a convex domain in $\mathbb{R}^n$. Let $u$ be an eigenfunction on $M$ with respect to the first Neumann eigenvalue $\mu_1(M)$. Suppose that $\sup u\geq-\inf u$.
	Let $\psi$ be the solution to the Sturm-Liouville problem
	\begin{align}
		\psi''+\frac{n-1}{x}\psi'+\psi=0,
	\end{align}
	with initial condition $\psi(0)=1,\psi'(0)=0$. Let	$C_n=-\Gamma(\frac{n}{2})(\frac{2}{x_m})^{\frac{n-2}{2}}J_{\frac{n-2}{2}}(x_m)$, where $x_m$ is the first minimum point of $t^{-\frac{n-2}{2}}J_{\frac{n-2}{2}}(t)$ for $t>0$, and $J_{\alpha}$ is Bessel function of order $\alpha$.
	
	Then, we have
	\begin{align}
		\frac{\sup u}{-\inf u}\leq\frac{\sup\psi}{-\inf\psi}=\frac{\sup_{t>0} t^{-\frac{n-2}{2}}J_{\frac{n-2}{2}}(t)}{-\inf_{t>0} t^{-\frac{n-2}{2}}J_{\frac{n-2}{2}}(t)}=\frac{1}{C_n}.
	\end{align}
\end{lemma}

\begin{lemma}\label{lem positive negative}
	Assume $u$ is the first Neumann eigenfunction of $\Omega$, $\Delta u= -\mu\cdot u$ on $\Omega\subseteq \mathbb{R}^n$. Assume $u(z_0)=0$ for some $z_0\in\bar{\Omega}$, then for any $\delta> 0$, there exist $p,q\in B_{z_0}(\delta)\cap \Omega$ such that $u(p)>0$ and $u(q)<0$.
\end{lemma}

\pf
By contradiction, we can further assume $u\leq 0$ on $B_{z_0}(\delta)\cap\Omega$ (otherwise we consider $-u\leq 0$ in the rest argument). Then we have 
\begin{align}
	\Delta u= -\mu\cdot u\geq 0,  \quad \quad \text{on}\ B_{z_0}(r)\cap \Omega. \nonumber 
\end{align}

If $u(z_1)=0$ for some $z_1\in B_{z_0}(\delta)\cap\Omega$, by strong maximum principle, we get $u$ is constant on $B_{z_0}(\delta)\cap\Omega$. By unique continuation property \cite[Theorem 1.1]{GL}, $u$ is constant in $\Omega$, which is the contradiction. 

Hence $u<0$ on $B_{z_0}(\delta)\cap\Omega$, we get $z_0\in \partial \Omega$ and 
\begin{align}
	u(z_0)>u(z),\,\forall z\in B_{z_0}(\delta)\cap\Omega. \nonumber 
\end{align}
By Hopf lemma, we get $\frac{\partial u}{\partial \vec{n}}>0$ at $z_0$, where $\vec{n}$ is the outward unit normal vector at $z_0$. It is the contradiction to the Neumann boundary condition of $u$.
\qed

\begin{lemma}\label{lem nodal set projection}
	Assume $\Omega\subset(0,d)\times\mathscr{P}_{v_0}(\Omega)$, and $A=(0,\mathbf{y}_0),\,B=(d,\mathbf{y}_1)\subset\overline{\Omega}$. Assume $\Delta u= -\mu\cdot u$ on $\Omega\subseteq \mathbb{R}^n$ and $u(A)=0$, define $\Sigma_p=\{x\in\mathbb{R}: \{x\}\times\mathbb{R}^{n-1}\cap u^{-1}(0)\neq\emptyset\}$, then $\Sigma_p=[0,a]$ for some $a\in(0,d]$.
\end{lemma}

\pf
Denote $\Sigma=u^{-1}(0)$. Since $u(A)=0$, by Lemma \ref{lem positive negative}, we get $u$ has positive and negative values in $B_A(\delta)\cap \Omega$ for any $\delta> 0$. So $u$ is zero at some point $q\in (B_A(\delta)\cap \Omega)$, so $\Sigma_p$ is not a single point. 

Next, we show $\Sigma_p$ is connected. Otherwise, there exists $0\leq a_1<a_2<a_3\leq d$, such that $a_1,a_3\in\Sigma_p$ and $a_2\notin\Sigma_p$.

Note $\Omega_{a_2}\vcentcolon= \{(a_2, \mathbf{y})\in \Omega: y\in \mathbb{R}^{n- 1}\}$ is connected, and $u\big|_{\Omega_{a_2}}\neq0$, without loss of generality, we can assume 
\begin{align}
	u\big|_{\Omega_{a_2}}> 0.   \label{u a2 value positive} 
\end{align}

Choose $p_1=(a_1,\mathbf{y}_1)$ and $p_3=(a_3,\mathbf{y}_3)$ with $p_1,p_3\in\Sigma$. By Lemma \ref{lem positive negative}, we can choose points $q_1$ and $q_3$ in the neighborhood $U_1, U_3$ of $p_1,p_3$, such that $u(q_1), u(q_3)<0$. Since $a_1<a_2<a_3$, the neighborhood $U_1,U_3$ can be chosen such that 
\begin{align}
	U_k\cap\{(x,y)\in\mathbb{R}^n,x=a_2\}=\emptyset, \quad \quad k= 1, 3. \nonumber 
\end{align}

Assume the first coordinate of $q_1$ and $q_3$ is $b_1$ and $b_3$, we have $b_1<a_2<b_3$. 

We claim $q_1$ and $q_3$ are not in the same nodal domain. Otherwise, there is a curve $\gamma$ in this nodal domain connecting $q_1$ and $q_3$, with
\begin{align}
	u\big|_{\gamma}< 0. \label{u gamma negative} 
\end{align}

The first coordinate of points in $\gamma$ starts from $b_1<a_2$ and end in $b_3>a_2$, so 
\begin{align}
	\gamma\cap \Omega_{a_2}\neq \emptyset. \label{gamma omega a2 non-empty inter} 
\end{align}
By (\ref{u a2 value positive}), (\ref{u gamma negative}) and (\ref{gamma omega a2 non-empty inter}), there is the contradiction. So $q_1$ and $q_3$ are in two different nodal domains. 

Since there is at least one nodal domain where $u> 0$, we obtain at least three nodal domains, which contradicts Courant-Cheng's Nodal domain Theorem.

Since $u(A)=0$, we get $\Sigma_p$ is a segment with an end point equal to $0$. The nodal set is a closed subset of $\bar{\Omega}$, so $\Sigma_p$ is a closed connected subset of $[0,d]$, and the conclusion follows. 
\qed

\begin{theorem}\label{thm higher thin domain has simple first eigenvalue}
	Let $\Omega\subset\mathbb{R}^n$ be a convex domain with $C^1$ smooth boundary, and width $\mathcal{W}(\Omega)=\epsilon$. If the first Neumann eigenvalue $\mu$ of $\Omega$ is not simple, let 	$C_n=-\Gamma(\frac{n}{2})(\frac{2}{x_m})^{\frac{n-2}{2}}J_{\frac{n-2}{2}}(x_m)$, $J_\alpha$ is Bessel function of order $\alpha$, and $x_m$ is the first minimum point of $t^{-\frac{n-2}{2}}J_{\frac{n-2}{2}}(t)$ for $t>0$, then we have
	\begin{align}
		\mu\geq\frac{C_n}{2\epsilon^2}.
	\end{align}
\end{theorem}

\pf
Assume $\Omega\subset(0,d)\times\mathscr{P}_{v_0}(\Omega)$, and $A=(0,\mathbf{y}_0),\,B=(d,\mathbf{y}_1)\subset\overline{\Omega}$. Because $\mu$ is not simple eigenvalue, by linear combination of two linear independent eigenfunctions of $\mu$, there exists an eigenfunction $u$ corresponding to $\mu$ such that $u(A)=0$, and $\sup u=1,\,\inf u=-k,\,0<k\leq1$.

Assume $P,Q\in\overline{\Omega}$ such that 
\begin{align}
	u(P)=\min u=-k, \quad u(Q)=\max u=1, \quad P=(p_1,\mathbf{p}),\quad Q=(q_1,\mathbf{q}). \nonumber 
\end{align}

Since $u(P)<0$ and $u(Q)>0$, and $\Omega$ is convex, we get there exists a point $R\in\Omega$ on segment $\overline{PQ}$, such that 
\begin{align}
	u(R)=0, \quad \quad R=(r_1,\mathbf{r}), \quad \quad \min\{p_1,q_1\}\leq r_1\leq\max\{p_1,q_1\}. \nonumber 
\end{align}

From Lemma \ref{lem nodal set projection}, we know that $r_1\in\Sigma_p=[0,a]$. 

If $p_1\leq q_1$, we get $p_1\leq r_1\leq a$ and $p_1\in\Sigma_p$, so $u(p_1,\mathbf{y})=0$ for some $(p_1,\mathbf{y})\in\Omega$. By Lemma \ref{us0 vertical control}, we get $u(p_1,\mathbf{y})\leq-k+2\mu\epsilon^2$, so we get $\mu\geq\frac{k}{2\epsilon^2}$.

If $q_1<p_1$, similar arguments show that $\mu\geq\frac{1}{2\epsilon^2}$.

In both cases, we get

\begin{align}\label{lower bound mu 1}
	\mu\geq\frac{k}{2\epsilon^2}.
\end{align}

By Lemma \ref{max min value mfd}, we get $k\geq=\frac{-\inf t^{-\frac{n-2}{2}}J_{\frac{n-2}{2}}(t)}{\max t^{-\frac{n-2}{2}}J_{\frac{n-2}{2}}(t)}=-\Gamma(\frac{n}{2})(\frac{2}{x_m})^{\frac{n-2}{2}}J_{\frac{n-2}{2}}(x_m)=C_n$, where $x_m$ is the first minimum point of $t^{-\frac{n-2}{2}}J_{\frac{n-2}{2}}(t)$. 

So, we get $\mu\geq\frac{C_n}{2\epsilon^2}$.

\qed

From Theorem \ref{thm higher thin domain has simple first eigenvalue} and Theorem \ref{upper bound of eigen of N} we get the following result,

\begin{theorem}\label{ndim eigen simple}
	Let $\Omega\subset\mathbb{R}^n$ be a convex domain with $C^1$ smooth boundary, and diameter $\mathcal{D}(\Omega)$, width $\mathcal(\Omega)$, if $\frac{\mathcal{W}(\Omega)}{\mathcal{D}(\Omega)}<\frac{\sqrt{C_n}}{\sqrt{8}j_{\frac{n-2}{2},1}}$, where $C_n=-\Gamma(\frac{n}{2})(\frac{2}{x_m})^{\frac{n-2}{2}}J_{\frac{n-2}{2}}(x_m)$, and $x_m$ is the first minimum point of $t^{-\frac{n-2}{2}}J_{\frac{n-2}{2}}(t)$. Then the first Neumann eigenvalue of $\Omega$ is simple.
\end{theorem}


\section{Simplicity of eigenvalue on manifolds}

\begin{lemma}\label{two part}
Let $M$ be a compact complete $n$-dim Riemannian manifold with $Rc\geq 0$. Let $A$, $B$ be two points in $M$, such that $d(A,B)=\mathrm{diam}\,M=d$, let $l$ be a segment in $M$ connecting $A$, $B$, such that 
\begin{align}
l:[0,d]\to M,\quad l(0)=A,\quad l(d)=B. \nonumber
\end{align}
Suppose that $M$ is in the $\epsilon$ neighborhood of $l$, i.e. for any $p\in M$, we have $d(p,l)<\epsilon$. Then for any $\alpha\in[0,d]$, we have 
\begin{align}
M\setminus B_{2\epsilon}(l(\alpha))\subset B_\alpha(A)\cup B_{d-\alpha}(B)
\end{align}
\end{lemma}

\pf
For any $z\in M\setminus B_{2\epsilon}(l(\alpha))$, since $d(z,l)<\epsilon$, there exists $l(\beta)\in l$, such that $d(z,l(\beta))<\epsilon$, since $d(z,l(\alpha))\geq 2\epsilon$, we get $\beta<\alpha-\epsilon$ or $\beta>\alpha+\epsilon$. In either case, we get $d(A,l(\beta))<\alpha-\epsilon$ or $d(B,l(\beta))<d-\alpha-\epsilon$. So $d(A,z)<\alpha$ or $d(B,z)<d-\alpha$.
\qed

\begin{theorem}\label{manifold simple eigenvalue}
Let $M$ be a compact complete $n$-dim Riemannian manifold with $Rc\geq 0$. Let $A$, $B$ be two points in $M$, such that $d(A,B)=\mathrm{diam}\,M=d$, let $l$ be a segment in $M$ connecting $A$, $B$, such that 
\begin{align}
l:[0,d]\to M,\quad l(0)=A,\quad l(d)=B. \nonumber
\end{align}
Suppose that $M$ is in the $\epsilon$ neighborhood of $l$, i.e. for any $p\in M$, we have $d(p,l)<\epsilon$. If the first eigenvalue $\mu$ of $M$ is not simple, then we have $\mu$ satisfies 
\begin{align}
\mu\geq\frac{C_n}{18\epsilon^2}.
\end{align}
where $C_n=-\Gamma(\frac{n}{2})(\frac{2}{x_m})^{\frac{n-2}{2}}J_{\frac{n-2}{2}}(x_m)$, $J_\alpha$ is Bessel function of order $\alpha$, and $x_m$ is the first minimum point of $t^{-\frac{n-2}{2}}J_{\frac{n-2}{2}}(t)$ for $t>0$,
\end{theorem}

\pf
Step(1): Since $\mu_1(M)$ is not simple, let $u_1$ and $u_2$ be two linear independent eigenfunction of $M$ with respect to $\mu_1(M)$. By linear combination of $u_1$ and $u_2$, we get an eigenfunction $u$ corresponding to $\mu_1(M)$, such that $u(A)=0$, and $\sup u=1, \,\inf u=-k,\,0<k\leq1$.

Let $v=\frac{u-\frac{1-k}{2}}{\frac{1+k}{2}}$ and $\theta=\sin^{-1}v\in[-\frac{\pi}{2},\frac{\pi}{2}]$. From \cite{ZY}, we have
\begin{align}
	|D\theta|\leq\sqrt{\mu}\sqrt{\frac{2}{1+k}}\leq\sqrt{2\mu}.
\end{align}
and we get
\begin{align}\label{Du est mfd}
	|Du|=\frac{1+k}{2}|Dv|=\frac{1+k}{2}\sqrt{1-v^2}|D\theta|\leq\sqrt{(u+k)(1-k)}\sqrt{2\mu}.
\end{align}
Let $u(P)=\min u=-k$ and $u(Q)=\max u=1$. We apply similar  argument as in Lemma \ref{us0 vertical control}, 

For $x\in B_{3\epsilon}(P)$, we have 
\begin{align}
	u(x)\leq u(P)+3\epsilon\sup|Du|\leq-k+3\epsilon\sqrt{2\mu}.
\end{align}
By (\ref{Du est mfd}), we get for $x\in B_{3\epsilon}(P)$, we have
\begin{align}
	|Du|(x)\leq \sqrt{2\mu}\sqrt{3\epsilon\sqrt{2\mu}}.
\end{align}
Repeat this process as in Lemma \ref{us0 vertical control}, we get for $x\in B_{3\epsilon}(P)$, we have
\begin{align}\label{u est mfd}
u(x)\leq -k+(3\epsilon\sqrt{2\mu})^{-1-\frac{1}{2}-\frac{1}{4}-\cdots}=-k+18\mu\epsilon^2.
\end{align}
Similarly, we get for $x\in B_{3\epsilon}(Q)$, we have 
\begin{align}
	u(x)\geq1-18\mu\epsilon^2.
\end{align}

By Lemma \ref{max min value mfd}, we get 
\begin{align}
	k\geq C_n,
\end{align}

If $18\mu\epsilon^2<C_n\leq k$, then we get $-k+18\mu\epsilon^2<0$. So we have $u<0$ on $B_{3\epsilon}(P)$, and $u>0$ on $B_{3\epsilon}(Q)$.

Step(2):
Let $\Sigma$ be the nodal set of $u$ and $\Sigma_0=\{z\in l,\,d(z,\Sigma)\leq 2\epsilon\}$. For any $\alpha\in[0,d]$, if $l(\alpha)\notin\Sigma_0$, then $u\neq 0$ on $B_{2\epsilon}(l(\alpha))$, assume $u>0$ on $B_{2\epsilon}(l(\alpha))$. 
Then, by lemma \ref{two part} we have 
\begin{align}
\{u<0\}\subset M\setminus B_{2\epsilon}(l(\alpha))\subset [B_\alpha(A)\cup B_{d-\alpha}(B)]. \nonumber
\end{align}
By Courant-Cheng's Nodal domain theorem \cite[Section 1]{Cheng}, we get $u$ has two nodal domain, $\{u<0\}$ and $\{u>0\}$. So $\{u<0\}$ is connected.

For any neighborhood $U$ of point $A$ in $\Sigma$, if $u\geq0$ on $U$, then we have 
\begin{align}
\Delta u=-\mu u\leq0   \nonumber
\end{align}
Since $u(A)=0$, we have $u(A)=\min_{z\in U}u(z)$. By strong maximum principle, we get $u$ is constant in $U$. So $u$ is constant function, contradiction. So $\{u<0\}\cap U\neq\emptyset$, so we get
\begin{align}
\{u<0\}\cap B_\alpha(A)\neq\emptyset.\nonumber
\end{align}

Since $B_\alpha(A)$ and $B_{d-\alpha}(B)$ is disconnected, we get $\{u<0\}\subset B_\alpha(A)$. and we get $\Sigma\subset \overline{B_\alpha(A)}$,

Step(3):
Choose $P_1=l(p)$ and $Q_1=l(q)$, such that $d(P,l)=d(P,P_1)<\epsilon$ and $d(Q,l)=d(Q,Q_1)<\epsilon$, Since $u<0$ on $B_{3\epsilon}(P)$, and $u>0$ on $B_{3\epsilon}(Q)$. We get $u<0$ on $B_{2\epsilon}(P_1)$ and $u>0$ on $B_{2\epsilon}(Q_1)$, so $P_1=l(p),Q_1=l(q)\notin\Sigma_0$. By Step(2), we get 
\begin{align}
	\Sigma\subset \overline{B_p(A)},\quad \Sigma\subset \overline{B_q(A)}.
\end{align}

However, let $\gamma$ be a curve consists of three segments in $M$, $\gamma=\overline{PP_1}\cup\overline{P_1Q_1}\cup\overline{Q_1Q}$. Since $u(P)<0$ and $u(Q)>0$, there is a point $x_0$ on $\gamma$ such that $u(x_0)=0$. Since $u<0$  on $B_{3\epsilon}(P)$, and $u>0$ on $B_{3\epsilon}(Q)$. we get $x_0\in\Sigma$ is on inner point of segment $\overline{P_1Q_1}$. So $d(x_0,A)>\min\{p,q\}$, contradict to $\Sigma\subset \overline{B_p(A)},\quad \Sigma\subset \overline{B_q(A)}$. 
So we get
\begin{align}
	18\mu\epsilon^2\geq C_n.
\end{align}
So 
\begin{align}
	\mu\geq\frac{C_n}{18\epsilon^2}.
\end{align}

\qed

\part{The Comparison of Neumann eigenvalues}

\section{The bounds of eigenvalues}

\begin{lemma}\label{lem put convex domain into Rn}
	{Let $\Omega\subset\mathbb{R}^n$ be a convex domain with diameter $\mathcal{D}(\Omega)=1$ and width $\mathcal{W}(\Omega)=\epsilon$. Suppose that $(0,\mathbf{0}), (1,\mathbf{0})\in \overline\Omega$, then 
		\begin{align}
			\Omega\subset \Big([0,1]\times B_{\mathbf{0}}(\epsilon)\Big); \nonumber 
		\end{align}
		where $B_{\mathbf{0}}(\epsilon)$ is the $(n-1)$-dimensional ball of radius $\epsilon$ centred at $\mathbf{0}$ in $\mathbb{R}^{n-1}$.
	}
\end{lemma}

\pf
{Denote the two points, $v_1=(0,\mathbf{0})$, and $v_2=(1,\mathbf{0})$.  Suppose there exists a point $v_3=(x,\mathbf{y})\in\Omega$, such that $|\mathbf{y}|=\epsilon_0>\epsilon$. 
	
	For any $v\in\mathbb{S}^{n-1}$, without loss of generality, we can assume $v\cdot v_i\leq v\cdot v_j\leq v\cdot v_k$, where $i,j,k$ is some permutation of $1,2,3$. Then there exists $\lambda\in[0,1]$ such that $v\cdot v_j=\lambda v\cdot v_i+(1-\lambda)v\cdot v_k=v\cdot(\lambda v_i+(1-\lambda)v_k)$. Note $v_0\vcentcolon=\lambda v_i+(1-\lambda)v_k$ is a point on segment $v_iv_k$. Since $\Omega$ is convex, we know that $v_0\in\Omega$.
	
	The distance of $v_3$ to line $v_1v_2$ is $|\mathbf{y}|$, so the area of the triangle formed by $v_1,v_2,v_3$ equals to 
	\begin{align}
		\mathscr{A}=\frac{1}{2}|v_1-v_2|\cdot|\mathbf{y}|=\frac12|\mathbf{y}|=\frac{\epsilon_0}{2}. \nonumber 
	\end{align}
	
	On the other hand, assume the angle between $\vec{v_0v_k}$ and $\vec{v_0v_j}$ is $\theta\in[0,2\pi)$. Since $\mathcal{D}(\Omega)=1$, we have $|v_i-v_k|\leq 1$. Then
	\begin{align}
		\epsilon_0= 2\mathscr{A}=|v_i-v_k|\cdot|v_j-v_0|\cdot |\sin \theta|\leq |v_j-v_0|. \nonumber 
	\end{align}
	
	Let $x_j=v_j-(v_j\cdot v)\cdot v$, and $x_0=v_0-(v_0\cdot v)\cdot v$, we have $x_j, x_0\in \mathscr{P}_{v}(\Omega)$. By $v\cdot v_j=v\cdot v_0$, we get $x_j-x_0=v_j-v_0$. Now we have
	\begin{align}
		\mathcal{D}(\mathscr{P}_{v}(\Omega))\geq |x_j-x_0|=|v_j-v_0|\geq\epsilon_0. \nonumber
	\end{align}
	
	Therefore $\mathcal{W}(\Omega)= \min_{v\in \mathbb{S}^{n- 1}}  \mathcal{D}(\mathscr{P}_{v}(\Omega))\geq \epsilon_0>\epsilon$, which contradicts $\mathcal{W}(\Omega)=\epsilon$.
	
	Now for any point $P=(x,\mathbf{y})\in\Omega$, we have $|\mathbf{y}|\leq\epsilon$. Since $|P-v_1|,|P-v_2|\leq\mathcal{D}(\Omega)=1$, we get $x\in[0,1]$. The conclusion follows.
}
\qed

\textbf{In the rest of this chapter, we assume $\Omega\subset\mathbb{R}^n$ is a convex domain with 
\begin{align}
&\mathcal{D}(\Omega)= 1, \quad\quad \quad (0,\mathbf{0}), (1,\mathbf{0})\in \overline\Omega, \quad \quad \mathcal{W}(\Omega)= \epsilon\in (0, \frac{1}{2}). \nonumber 
\end{align}}

From Lemma \ref{lem put convex domain into Rn}, we know 
\begin{align}
\Omega\subset \Big([0,1]\times B_{\mathbf{0}}(\epsilon)\Big); \nonumber 
\end{align}

Assume $\Omega_t= \{y\in \mathbb{R}^{n-1}: (t, y)\in \Omega\}$ and $\Omega=\{(t,y):\,t\in[0,1],\,y\in\Omega_t\}$. We have $\Omega_t\subset \{t\}\times B_{\mathbf{0}}(\epsilon)$.

Let $V= V(\Omega)=\int_{0}^1\int_{\Omega_x}1\,dt$ in the rest argument.

For $L^2(\Omega)$,  we use the norm as follows:
\begin{align}
|f|_{L^2(\Omega)}\vcentcolon = \Big(\frac{1}{V}\int_{\Omega}f(x, y)^2dxdy \Big)^{\frac{1}{2}}.  \nonumber 
\end{align}

We define
\begin{align}
H(x)\vcentcolon= \frac{\mathcal{H}^{n-1}(\Omega_x)}{\int_{0}^1 \mathcal{H}^{n-1}(\Omega_t)dt},  \nonumber 
\end{align}
where $\mathcal{H}^{n- 1}$ is $(n-1)$-dim Hausdorff measure.

We assume $d\mu(x)= H(x)dx$,  then $\mu$ is a Radon measure on $N\vcentcolon= [0,  1]$.  

For $L^2(N)$,  we use the norm as follows:
\begin{align}
|h|_{L^2(N)}\vcentcolon = \Big(\int_{0}^1 h(x)^2d\mu(x) \Big)^{\frac{1}{2}}.  \nonumber 
\end{align} 

The min-max definition of the Neumann eigenvalues $\Omega,  N$ yields: for $k\geq 1$,  
\begin{align}
\mu_{k- 1}(\Omega)=& \inf_{V_k\subseteq C^\infty(\Omega)}\sup_{f\in V_k-\{0\}}\frac{\int_{\Omega}|Df|^2}{\int_{\Omega} f^2},   \label{lambda k Omega i} \\
\mu_{k- 1}(N)=& \inf_{V_k\subseteq C^\infty(N)}\sup_{h\in V_k-\{0\}}\frac{\int_{N}|Dh|^2}{\int_{N} h^2}.   \label{lambda k N} 
\end{align}
where $V_k$ is any $k$-dim subspace.

\begin{definition}\label{def lambda for k-dim subspace}
{For $k$-dim subspace $\Lambda^k$, we define
\begin{align}
\lambda(\Lambda^k)= \sup_{f\in \Lambda^k- \{0\}}\frac{\int |Df|^2}{\int f^2}. \nonumber 
\end{align}
}
\end{definition}

We define $\Psi: L^2(N)\rightarrow L^2(\Omega)$ as follows: for any $h\in L^2(N)$,
\begin{align}
\Psi(h)(x, y)\vcentcolon= h(x), \quad \quad \quad \forall x\in [0, 1]\ \text{and}\ (x, y)\in \Omega. \nonumber 
\end{align}

\begin{prop}\label{prop upper bound of mu}
{For any $k\in \mathbb{Z}^+$, we have $\displaystyle \mu_{k- 1}(\Omega)\leq \mu_{k- 1}(N)$. 
}
\end{prop}

\pf
{Note for any $\delta> 0$,  from (\ref{lambda k N}),  we can find $k$-dim subspace $\Lambda^k\subseteq C^\infty(N)$ such that 
\begin{align}
\lambda(\Lambda^k)\leq \mu_{k- 1}(N)+ \delta.  \label{1st eigenvalue eq}
\end{align}

Put $\tilde{\Lambda}^k= \Psi(\Lambda^k)$.  Note $\Psi$ is linear and $\mathrm{Ker}(\Psi)=\{0\}$, so $\tilde{\Lambda}^k\subseteq C^\infty(\Omega)$ is a $k$-dim subspace.  By (\ref{lambda k Omega i}) we get
\begin{align}
\mu_{k- 1}(\Omega)\leq \lambda(\tilde{\Lambda}^k).  \label{second eigenvalue eq}
\end{align}

Assume there is $h_0\in \Lambda^k$ such that 
\begin{align}
\lambda(\tilde{\Lambda}^k)= \frac{\int_\Omega |D(\Psi(h_0))|^2}{\int_\Omega |\Psi(h_0)|^2} . \nonumber 
\end{align}

Then note $\displaystyle \frac{\int_\Omega |D(\Psi(h_0))|^2}{\int_\Omega |\Psi(h_0)|^2}= \frac{\int_N |D(h_0)|^2}{\int_N |h_0|^2}$, using (\ref{1st eigenvalue eq}) and (\ref{second eigenvalue eq}); we get 
\begin{align}
\mu_{k- 1}(\Omega)\leq \lambda(\tilde{\Lambda}^k)\leq \frac{\int_N |D(h_0)|^2}{\int_N |h_0|^2} \leq \lambda(\Lambda^k)\leq \mu_{k- 1}(N)+ \delta. \nonumber 
\end{align}

Let $\delta\rightarrow 0$ in the above, the conclusion follows.
}

We define $\Phi: L^2(\Omega)\rightarrow L^2(N)$ as follows: for any $f\in L^2(\Omega)$, 
\begin{align}
\Phi(f)(x)\vcentcolon= f(x, 0), \quad \quad \quad \forall x\in N. \nonumber 
\end{align}
It is obvious that $\Phi\circ \Psi= \mathrm{Id}: L^2(N)\rightarrow L^2(N)$.

\begin{lemma}\label{lem varing domain integal ests}
{For any $f\in C^\infty(\Omega)$, then 
\begin{align}
\Big||f|_{L^2(\Omega)}^2- |\Phi(f)|_{L^2(N)}^2\Big|\leq  \epsilon\cdot \Big\{ |f|_{L^2(\Omega)}^2+ |Df|_{L^2(\Omega)}^2\Big\}. \nonumber 
\end{align}
}
\end{lemma}

\pf
We have 

\begin{align}
&\Big||f|_{L^2(\Omega)}^2-|\Phi(f)|_{L_2(N)}^2\Big| \nonumber \\
=&\frac{1}{V}\Big|\int_0^1dx\int_{\Omega_x}(f(x,y)^2-f(x,0)^2)\,dy\Big|   \nonumber \\
=&\frac{1}{V}\Big|\int_0^1dx\int_{\Omega_x}dy\int_0^y\partial_t(f(x,t)^2)\,dt\Big|   \nonumber \\
\leq&\frac{1}{V}\Big|\int_0^1dx\int_{\Omega_x}dy\int_0^y2\cdot|f(x,t)|\cdot|Df(x,t)|\,dt\Big|   \nonumber \\
\leq&\frac{1}{V}\Big|\int_0^1dx\int_{\Omega_x}\epsilon|2f(x,y)|\cdot |Df(x,y)|dy \Big|  \nonumber \\
\leq&\frac{\epsilon}{V}\Big|\int_0^1dx\int_{\Omega_x}f(x,y)^2+|Df(x,y)|^2\,dt\Big|   \nonumber \\
=&\epsilon(|f|_{L^2(\Omega)}^2+|Df|_{L^2(\Omega)}^2)
\end{align}
\qed

The following lemma is well-known from the convexity property of $\Omega$. 
\begin{lemma}\label{lem boundary deriv est of eigenfunc}
{For any $w\in C^\infty(\bar{\Omega})$ with $\frac{\partial w}{\partial \vec{n}}\big|_{\partial \Omega}= 0$, we have $\displaystyle \Big(D_{\vec{n}}D_{\vec{t}}w\cdot D_{\vec{t}}w\Big)\big|_{\partial \Omega}\leq 0$, where $\vec{t}$ is unit tangent vector of $\partial \Omega$.
}
\end{lemma}\qed

\begin{lemma}\label{lem Hessian est by Lap on convex for Neumann func}
{For any $u\in C^\infty(\Omega)$ with $\frac{\partial u}{\partial \vec{n}}\big|_{\partial\Omega}= 0$,  we have 
\begin{align}
\int_\Omega |D^2 u|^2\leq \int_\Omega |\Delta u|^2.  \nonumber 
\end{align}
}
\end{lemma}

\pf
{Let $e_1,\cdots,e_{n-1}$ be orthonormal frame on $\partial\Omega$.
By Lemma \ref{lem boundary deriv est of eigenfunc},  we have $D_iu\cdot D_i\frac{\partial u}{\partial n}\leq 0$, for $i=1,\cdots,n-1$.
\begin{align}
\frac{1}{2}\int_{\Omega}\Delta(|Du|^2)= \frac{1}{2}\int_{\partial\Omega}\frac{\partial}{\partial\vec{n}}|Du|^2= \int_{\partial\Omega}Du\cdot D\frac{\partial u}{\partial n}=\int_{\partial\Omega}\sum_{i=1}^{n-1} D_iu\cdot D_i\frac{\partial u}{\partial n}\leq 0.  \nonumber 
\end{align}

On the other hand, 
\begin{align}
\frac{1}{2}\int_{\Omega}\Delta(|Du|^2)= \int_{\Omega}|D^2u|^2+ D(\Delta u)\cdot Du= \int_{\Omega}|D^2u|^2- |\Delta u|^2.  \nonumber 
\end{align}
}
\qed

\begin{lemma}\label{lem varing domain integal ests for grad}
{For any $f\in C^\infty(\Omega)$, then 
\begin{align}
|Df|_{L^2(\Omega)}^2- |D(\Phi(f))|_{L^2(N)}^2\geq  -\epsilon\cdot \Big\{ |Df|_{L^2(\Omega)}^2+ |D^2 f|_{L^2(\Omega)}^2\Big\}. \nonumber 
\end{align}
}
\end{lemma}

\pf
{By Lemma \ref{lem varing domain integal ests},  we have
\begin{align}
&|Df|_{L^2(\Omega)}^2- |D(\Phi(f))|_{L^2(N)}^2\geq \frac{1}{V}\Big\{\int_0^1 dx\int_{\Omega_x} |Df(x, y)|^2 dy- \int_0^1 f_x^2(x, 0)\mathcal{H}^{n- 1}(\Omega_x)dx \Big\} \nonumber \\
&\geq \frac{1}{V}\Big\{\int_0^1 dx\int_{\Omega_x} |f_x(x, y)|^2 dy- f_x^2(x, 0)dy \Big\} \nonumber \\
&\geq -\epsilon\cdot \Big\{ |Df|_{L^2(\Omega)}^2+ |D^2 f|_{L^2(\Omega)}^2\Big\}.  \nonumber
\end{align}
}
\qed

\begin{prop}\label{prop inf limit is more}
{For any $k\geq 0$, we have 
\begin{align}
\mu_k(\Omega )\geq \Big\{1-2\epsilon(1+\mu_k(N))\Big\}\cdot\mu_k(N).\nonumber 
\end{align}
}
\end{prop}

\pf
{When $k= 0$, the conclusion is trivial. In the rest argument $k\geq 1$ is fixed. Note $\displaystyle \mu_k(\Omega)\leq \mu_k(N)$ by Proposition \ref{prop upper bound of mu}. 

Let $\{\varphi_{j}\}_{j= 1}^\infty$ be orthonormal normalized eigenfunctions of $\Omega$ and the corresponding eigenvalues $\mu_j(\Omega)$ is increasing in $j$. Put $\Lambda\vcentcolon= \mathrm{span}\{\varphi_{j}\}_{j= 1}^k$, 

From $f\in \Lambda$ and $|f|_{L^2(\Omega)}= 1$, we know that $\frac{\partial f}{\partial \vec{n}}\big|_{\partial\Omega}= 0$. By Lemma \ref{lem Hessian est by Lap on convex for Neumann func} and $\Delta \varphi_{j}= -\mu_j(\Omega)\varphi_{ j}$, we get 
\begin{align}\label{Df control}
&|D^2 f|_{L^2(\Omega)}^2\leq |\Delta f|_{L^2(\Omega)}^2\leq \mu_k(\Omega)^2|f|^2_{L^2(\Omega)}= \mu_k(\Omega)^2. \nonumber \\
&|Df|_{L^2(\Omega)}^2\leq \mu_k(\Omega)|f|^2_{L^2(\Omega)}\leq \mu_k(\Omega). 
\end{align}

We assume $1-2\epsilon(1+\mu_k(N))>0$, otherwise the right hand side of the proposition is negative. By Lemma \ref{lem varing domain integal ests},  Proposition \ref{prop upper bound of mu} and (\ref{Df control}) we get $\mathrm{Ker}(\Phi)=\{0\}$, so $\Phi(\Lambda)$ is a $k$-dimensional subspace of $C^\infty(N)$.

Choose $f\in \Lambda$ with $|f|_{L^2(\Omega)}= 1$ such that
\begin{align}
\lambda(\Phi(\Lambda))= \frac{\int |D\Phi(f)|^2}{\int \Phi(f)^2}. \nonumber 
\end{align}
From Lemma \ref{lem varing domain integal ests for grad}, we have
\begin{align}
\mu_k(\Omega )= \lambda(\Lambda )\geq \frac{|Df|_{L^2(\Omega)}}{|f|_{L^2(\Omega)}}\geq \frac{|D(\Phi (f ))|_{L^2(N)}^2 -\epsilon \cdot \Big\{ |Df |_{L^2(\Omega )}^2+ |D^2 f |_{L^2(\Omega )}^2\Big\}}{|f|_{L^2(\Omega)}}. \label{intermediate ineq-1}
\end{align}

From Lemma \ref{lem varing domain integal ests}, we get 
\begin{align}
|\Phi (f )|_{L^2(N)}^2\geq |f|_{L^2(\Omega)}^2- \epsilon\cdot \Big\{ |f |_{L^2(\Omega )}^2+ |Df |_{L^2(\Omega )}^2\Big\}\geq 1- \epsilon(1+ \mu_k(\Omega)). \nonumber 
\end{align}

Plugging the above inequalities into (\ref{intermediate ineq-1}), we get
\begin{align}
\mu_k(\Omega )& \geq \frac{|D(\Phi (f ))|_{L^2(N)}^2 - \epsilon \cdot (\mu_k(\Omega)^2+ \mu_k(\Omega))}{|f |_{L^2(\Omega )}^2}\nonumber \\
&= |D(\Phi (f))|_{L^2(N)}^2 - \epsilon \cdot (\mu_k(\Omega)^2+ \mu_k(\Omega)) \nonumber \\
&\geq \frac{|D(\Phi (f ))|_{L^2(N)}^2}{|\Phi (f )|_{L^2(N)}^2}\cdot |\Phi (f )|_{L^2(N)}^2 - \epsilon \cdot (\mu_k(\Omega)^2+ \mu_k(\Omega))\nonumber \\
&\geq \Big\{1- \epsilon(1+ \mu_k(\Omega))\Big\}\cdot \lambda(\Phi (\Lambda) )- \epsilon \cdot (\mu_k(\Omega)^2+ \mu_k(\Omega))\nonumber \\
&\geq \Big\{1- \epsilon(1+ \mu_k(\Omega))\Big\}\cdot \mu_k(N)- \epsilon \cdot (\mu_k(\Omega)^2+ \mu_k(\Omega))\nonumber \\
&\geq \Big\{1-2\epsilon(1+\mu_k(N))\Big\}\cdot\mu_k(N). \nonumber
\end{align}
}
\qed

\begin{theorem}\label{thm Omegai and N comparsion}
{Let $\Omega\subset\mathbb{R}^n$ be a convex domain with diameter $\mathcal{D}(\Omega)=1$ and width $\mathcal{W}(\Omega)=\epsilon$. Then, we have 
\begin{align}
\mu_k(N)\geq \mu_k(\Omega )\geq \Big\{1- 2\epsilon(1+ \mu_k(N))\Big\}\cdot \mu_k(N). \nonumber 
\end{align}
}
\end{theorem}

\pf
{It follows from Proposition \ref{prop upper bound of mu} and Proposition \ref{prop inf limit is more}. 
}
\qed

The positive zeros of the Bessel function $J_v$ of order $v$ are in increasing order 
\begin{align}
j_{v, 1}< j_{v, 2}<\cdots . \nonumber  
\end{align}

We recall the following result from \cite{Kr}.
\begin{prop}\label{prop upper bound of eigen of N}
{Let $\Omega\subset\mathbb{R}^n$ be a convex domain with diameter $\mathcal{D}(\Omega)=1$,let $k\in\mathbb{N}^+$, then we have
\begin{enumerate}
\item[(a)]. if $n= 2$, then $\mu_k(N)\leq 4(j_{0, 1}+ \frac{k- 1}{2}\pi)^2$;
\item[(b)]. if $n> 2$, and $k$ is odd, then $\mu_k(N)\leq 4j_{\frac{n-2}{2}, \frac{k+ 1}{2}}^2$;
\item[(c)]. if $n> 2$, and $k$ is even, then $\mu_k(N)\leq (j_{\frac{n-2}{2}, \frac{k}{2}}+ j_{\frac{n-2}{2}, \frac{k+ 2}{2}})^2$;
\end{enumerate}
}
\end{prop}

\pf
{See \cite[Proposition $2$]{Kr}.
}
\qed

\begin{cor}\label{cor 1st eignvalue est}
{Let $\Omega\subset\mathbb{R}^n$ be a convex domain with diameter $\mathcal{D}(\Omega)=1$ and width $\mathcal{W}(\Omega)=\epsilon$. There is $C(n,k)> 0$, such that 
\begin{align}
\mu_k(N)\geq \mu_k(\Omega )\geq \Big\{1- C(n,k)\cdot \epsilon)\Big\}\cdot \mu_k(N). \nonumber 
\end{align}
}
\end{cor}

\begin{remark}\label{rem another proof of Jerison's result}
{This proposition is due to \cite{Jerison} for $n= 2,\, k=1$, and for $n\geq 2, k= 1$ see \cite{CJK}. The proof we provide here is different from theirs. 
}
\end{remark}

\pf
{From Theorem \ref{thm Omegai and N comparsion}, we get
\begin{align}
\mu_k(\Omega)\geq (1- 2(1+ \mu_k(N))\cdot \epsilon)\cdot \mu_k(N). \nonumber 
\end{align}
From Proposition \ref{prop upper bound of eigen of N}, we get $\mu_k(N)$ is bounded by constant depending on $n,k$.
The conclusion follows.
}
\qed

\section{Vertical derivative estimate}
In this section, we assume $\Omega$ is a convex domain in $\mathbb{R}^2$, with diameter $\mathcal{D}(\Omega)=1$, and width $\mathcal{W}(\Omega)=\epsilon$ for some $\epsilon<\frac{1}{40}$. By Lemma \ref{lem put convex domain into Rn}, we further assume $\Omega\subseteq(0,1)\times(-\epsilon,\epsilon)\subseteq \mathbb{R}^2$, and $(0,0),(1,0)\in\overline{\Omega}$.

By \cite{Kr}, we have the first Neumann eigenvalue $\mu_1$ of $\Omega$ satisfies $\mu_1\leq 4j_{0,1}^2<25$ is bounded. Let $u$ be the first Neumann eigenfunction with respect to $\mu_1$, such that $\sup u=1$ and $\inf u=-k$, for $k\in(0,1]$, by Lemma \ref{max min value mfd}\cite{Kr2}, we get $k\geq \frac{1}{j_{0,1}}>\frac{2}{5}$.

\begin{lemma}\label{end point u Du control}
Let $u$ be the first Neumann eigenfunction of $\Omega$, with respect to eigenvalue $\mu_1$, assume $\sup u=1,\,\inf u=-k$, such that $0<k\leq1$, then
\begin{align}
\sup_{x\in[0,\epsilon]\cup[1-\epsilon,1]}|Du|(x,y)\leq 6\mu_1\epsilon,
\end{align}
we also have one of the following holds,
\begin{enumerate}
\item[(1)]. $\sup_{x\in[0,\epsilon]}u(x,y)\leq -k+18\mu_1\epsilon^2,\quad \inf_{x\in[1-\epsilon,1]}u(x,y)\geq1-18\mu_1\epsilon^2$ \\
\item[(2)]. $\inf_{x\in[0,\epsilon]}u(x,y)\geq 1-18\mu_1\epsilon^2,\quad \sup_{x\in[1-\epsilon,1]}u(x,y)\leq -k+18\mu_1\epsilon^2$.
\end{enumerate}
\end{lemma}

\pf
Let $v=\frac{u-\frac{1-k}{2}}{\frac{1+k}{2}}$ and $\theta=\sin^{-1}v\in[-\frac{\pi}{2},\frac{\pi}{2}]$. From \cite{ZY}, we have
\begin{align}
	|D\theta|\leq\sqrt{\mu_1}\sqrt{\frac{2}{1+k}}\leq\sqrt{2\mu_1}.
\end{align}
So, we get
\begin{align}\label{Du control}
|Du|=\frac{1+k}{2}|Dv|=\frac{1+k}{2}\sqrt{1-v^2}|D\theta|\leq\sqrt{(u+k)(1-k)}\sqrt{2\mu_1}.
\end{align}

Let $s_0=\max\{x\in[0,1],\min_y u(x,y)=-k\}$, then we get for any $(s_0,y)\in\Omega$, we have
\begin{align}\label{us0 control}
u(s_0,y)\leq-k+\sup|Du|\cdot 2\epsilon\leq-k+2\epsilon\sqrt{2\mu_1}.
\end{align}
Plug in (\ref{us0 control}) in (\ref{Du control}), we get
\begin{align}\label{Dus0 control}
|Du|(s_0,y)\leq \sqrt{2\mu_1}\sqrt{u(s_0,y)+k}\leq\sqrt{2\mu_1}(2\epsilon\sqrt{2\mu_1})^\frac{1}{2}.
\end{align}
Plug in (\ref{Dus0 control}) in (\ref{us0 control}) we get
\begin{align}\label{us0 control 2}
u(s_0,y)\leq-k+2\epsilon\cdot\sqrt{2\mu_1}(2\epsilon\sqrt{2\mu_1})^\frac{1}{2}.
\end{align}
Repeat this process $m$ times, by induction we get
\begin{align}
|Du|(s_0,y)\leq\sqrt{2\mu_1}(2\epsilon\sqrt{2\mu_1})^{2^{-1}+2^{-2}+\cdots+2^{-m}}.
\end{align}
Let $m\to\infty$, we get
\begin{align}
|Du|(s_0,y)\leq 4\mu_1\epsilon,
\end{align}
and
\begin{align}
u(s_0,y)\leq -k+8\mu_1\epsilon^2.
\end{align}

So, we have $u<0$ on $\Omega\cap\{x=s_0\}$. 

Since $u$ is the first eigenfunction, by Courant-Cheng Nodal domain theorem\cite{Cheng}, we get $\{u>0\}$ is connected, so 
\begin{align}
\{u>0\}\cap\{x<s_0\}=\emptyset,\quad or \quad \{u>0\}\cap\{x>s_0\}=\emptyset.
\end{align}
Without loss of generality, we assume
\begin{align}
\{u>0\}\cap\{x<s_0\}=\emptyset.
\end{align}

If $s_0=0$, then we get 
\begin{align}\label{(-1,0) control}
u(0,0)\leq-k+8\mu_1\epsilon^2.
\end{align}

If $s_0\neq0$, then on $\Omega_1\vcentcolon=(\Omega\cap\{x<s_0\})$ we have 
\begin{align}
\Delta u=-\mu_1u\geq0.
\end{align}
Assume $u(x_0)=\max_{x\in\bar{\Omega_1}}u(x)$. By maximum principle we get $x_0\in\partial\Omega_1$. If $x_0\in\mathrm{Interior}(\partial\Omega_1\cap\partial\Omega)$, by strong maximum principle, we get $\frac{\partial u}{\partial\vec{n}}(x_0)>0$, contradiction to Neumann boundary condition of $u$. So we get $x_0\in\bar{\Omega}\cap\{x=s_0\}$, and we get
\begin{align}\label{Omega1 max}
\max_{x\in\bar{\Omega}_1}u(x)\leq-k+8\mu_1\epsilon^2.
\end{align}

If $s_0>\epsilon$, by \ref{Omega1 max} we get 
\begin{align}
\sup_{x\in[0,\epsilon]}u(x,y)\leq\sup_{\bar{\Omega}_1}u\leq-k+8\mu_1\epsilon^2,
\end{align}
and by (\ref{Du control}), we get
\begin{align}
\sup_{x\in[0,\epsilon]}|Du(x,y)|\leq 4\mu_1\epsilon.
\end{align}

If $s_0\leq\epsilon$, assume $u(s_0,y_0)=-k$, then for any $(p,q)\in\Omega,x\in[0,\epsilon]$, we have $d((p,q),(s_0,y_0))\leq \epsilon+2\epsilon=3\epsilon$. So
\begin{align}
u(p,q)\leq u(s_0,y_0)+\sup_{x\in[0,\epsilon]}|Du|(x,z)\cdot 3\epsilon=-k+3\epsilon\cdot \sup_{x\in[0,\epsilon]}|Du|(x,z).
\end{align}
By (\ref{Du control}), we get
\begin{align}
|Du|(p,q)\leq \sqrt{2\mu_1}\sqrt{3\epsilon\sup_{x\in[0,\epsilon]}|Du|(x,z)}.
\end{align}
Taking sup for $p\in[0,\epsilon]$ on left side of the above inequality, we get 
\begin{align}
\sup_{x\in[0,\epsilon]}|Du|(x,y)\leq6\mu_1\epsilon,
\end{align}
and 
\begin{align}
\sup_{x\in[0,\epsilon]}u(x,y)\leq-k+18\mu_1\epsilon^2.
\end{align}
The same argument applied to $[1-\epsilon,1]$, we get the other half of the inequality.
\qed

\begin{lemma}\label{interior Du control}
We have
\begin{align}
\sup_{(x,y)\in\Omega}|D_yu(x,y)|\leq 2\sup_{(x,y)\in\partial\Omega}|D_yu(x,y)|.
\end{align}
\end{lemma}
\pf
Let $w(x,y)=1+\cos(\frac{y}{\epsilon}),\,v(x,y)=\frac{D_yu(x,y)}{w(x,y)}$.

If there is $(x_0,y_0)\in\Omega$ such that $|v(x_0,y_0)|=\sup_{(x,y)\in\Omega}|v(x,y)|$, without loss of generality we assume $v(x_0,y_0)>0$. We have
\begin{align}
0=\nabla v(x_0,y_0)=\frac{\nabla u_y}{w}-u_y\frac{\nabla w}{w^2},
\end{align}
and we have
\begin{align}
0&\geq\Delta v(x_0,y_0)=-\mu_1\frac{u_y}{w}+u_y(\frac{-\Delta w}{w^2}+2\frac{w_y^2}{w^3})-2\nabla u_y\frac{\nabla w}{w^2}\nonumber\\
&=(-\mu_1-\frac{\Delta w}{w})\cdot v(x_0,y_0).
\end{align}
By definition of $w$, we get
\begin{align}
0\geq(\frac{1}{\epsilon^2}-\mu_1)v(x_0,y_0)>0,
\end{align}
contradiction. So $|v(x,y)|\leq\sup_{(x,y)\in\partial\Omega}|v(x,y)|$.

Note $y\in[-\epsilon,\epsilon]$, we have $w(x,y)\in[1,2]$. So
\begin{align}
|D_yu(x,y)|\leq2|v(x,y)|\leq 2\sup_{(x,y)\in\partial\Omega}|v(x,y)|\leq 2\sup_{(x,y)\in\partial\Omega}|D_yu(x,y)|.
\end{align}

\qed

We have the following control of derivative of $u$ in vertical direction from \cite{Zheng}.

\begin{lemma}\label{vertical derivative control}\cite[Lemma $15$]{Zheng}
We have
\begin{align}
|D_yu(x,y)|\leq 48\mu_1\epsilon.
\end{align}
\end{lemma}

\pf
Step(1). We assume the first case in Lemma \ref{end point u Du control} happens, we have $\sup_{x\in[0,\epsilon]}u(x,y)\leq -k+18\mu_1\epsilon^2$.  If $x\in[0,\epsilon]$ or $x\in[1-\epsilon,1]$, by Lemma \ref{end point u Du control}, we get 
\begin{align}
|D_yu|\leq |Du|\leq 6\mu_1\epsilon.
\end{align}

Choose $m\in\mathbb{N}$, such that $\epsilon\in[2^{-m-1},2^{-m}]$. For $n=0,1,\cdots,m+1$, let
\begin{align}
|u|_n=\sup_{x\in[2^{-n-1},2^{-n}]}u(x,y),\quad |Du|_n=\sup_{x\in[2^{-n-1},2^{-n}]}|Du(x,y)|.
\end{align}
By Lemma \ref{end point u Du control}, we get
\begin{align}\label{induction base}
|u|_{m+1}\leq -k+18\mu_1\epsilon^2\leq-k+18\mu_1\cdot 2^{-2m}, \quad |Du|_{m+1}\leq 6\mu_1\epsilon\leq 6\mu_12^{-m}.
\end{align}

Step(2). We prove by induction that 
\begin{align}\label{induction u Du}
|u|_{n}\leq-k+18\mu_1\cdot2^{-2n+2},\quad |Du|_n\leq6\mu_1\cdot2^{-n+1},\quad n=0,1\cdots,m+1,
\end{align}
where $n=m+1$ is proved in \eqref{induction base}.

Note
\begin{align}
|u|_{n-1}\leq |u|_{n}+(2^{-n}+3\epsilon)|Du|_{n-1}\leq|u|_{n}+7\cdot 2^{-n}|Du|_{n-1},
\end{align}
and
\begin{align}
|Du|_{n-1}\leq\sqrt{2\mu_1}\sqrt{k+|u|_{n-1}}\leq\sqrt{2\mu_1}\sqrt{7\cdot 2^{-n}|Du|_{n-1}+18\mu_1\cdot2^{-2n+2}}.
\end{align}
Thus $|Du|_{n-1}\leq 6\mu_1\cdot2^{-n+2}$, and
\begin{align}
|u|_{n-1}\leq-k+18\mu_1\cdot2^{-2n+2}+42\mu_1\cdot2^{-2n+2}\leq-k+18\mu_1\cdot2^{-2n+4}.
\end{align}
So \eqref{induction u Du} is proved.
For $x\in[2^{-n-1},2^{-n}]$, where $n\leq m+1$, we have 
\begin{align}
|Du|(x,y)\leq |Du|_n\leq6\mu_1\cdot 2^{-n+1}\leq24\mu_1 x.
\end{align}
So we get for $x\in[\epsilon,\frac12]$, we have
\begin{align}
|Du|(x,y)\leq24\mu_1x.
\end{align}
Similar argument holds for $x\in[\frac12,1-\epsilon]$, so for $x\in[\epsilon,1-\epsilon]$, we have
\begin{align}\label{Du control 2}
|Du|(x,y)\leq 24\mu_1\min\{x,(1-x)\},\quad\forall (x,y)\in\Omega.
\end{align}

Step(3). 
For $p=(x,y)\in\partial\Omega$, let $l$ be a line tangent to $\Omega$ at $p$. Since $\Omega$ is convex, we get $(0,0),(1,0)$ lie on the same side of $l$. Without loss of generality, we assume $y>0$, let the slop for $l$ be $s$, we get
\begin{align}
s\in[-\frac{y}{1-x},\frac{y}{x}]\subset[-\frac{\epsilon}{1-x},\frac{\epsilon}{x}].
\end{align}
Since $\partial_nu(x,y)=0$, we get $D_yu(x,y)=s\cdot D_xu(x,y)$, by \eqref{Du control 2}, we get
\begin{align}\label{boundary Du control}
|D_yu(x,y)|\leq\frac{\epsilon}{\min\{x,1-x\}}|Du(x,y)|\leq 24\mu_1\epsilon
\end{align}
From Lemma \ref{interior Du control}, we get 
\begin{align}
|D_yu(x,y)|\leq 48\mu_1\epsilon.
\end{align}
\qed

\section{The $L^\infty$ estimate of Neumann eigenfunction}
In this section, we prove the $L^\infty$ of Neumann eigenfunction can be bounded by the $L^2$ norm with Sobolev-Neumann inequality. 

Assume $\Omega\in\mathbb{R}^2$ is a convex open domain and $\mathcal{D}(\Omega)=1$.

Define the \textbf{Sobolev-Neumann constant $SN_2(\Omega)$ for $\Omega$} as follows:
\begin{align}
	SN_2(\Omega)=V(\Omega)^{\frac12}\cdot \sup_{f\in C^\infty(\Omega)}\frac{|f- \fint_\Omega f|_{L^2}}{ |Df|_{L^1}}.  \nonumber 
\end{align}

From \cite{DWZ}[Theorem 1.4], we have the following upper bound for the  $SN_2(\Omega)$ which depends only on the diameter of $\Omega$.

\begin{theorem}\label{SN2 upper bound}
	Assume $\Omega\in\mathbb{R}^2$ is a convex open domain. Then
	$$
	SN_2(\Omega)\leq 1600\mathcal{D}(\Omega).
	$$
\end{theorem}

\begin{lemma}\label{lem Sobolev improved deri ests}
	{For all $f\in W^{1, 2}(\Omega)$, we have
		\begin{align}
			V(\Omega)^{\frac12}|f|_4^2\leq |f|_2^2+ 4(SN_2(\Omega))^2V(\Omega)\cdot |Df|_2^2 .  \nonumber 
		\end{align}
	}
\end{lemma}

\pf
{Let $F= f^2$,  firstly note 
	\begin{align}
		|DF|_1^2= 4|fDf|_1^2\leq 4|f|_2^2\cdot |Df|_2^2.  \nonumber 
	\end{align}
	
	On the other hand we have
	\begin{align}
		|DF|_1^2\geq (SN_2(\Omega))^{-2}|F- \fint_\Omega F|_2^2= (SN_2(\Omega))^{-2}\cdot \Big\{|f|_4^4- V(\Omega)^{-1}|f|_2^4\Big\}.  \nonumber 
	\end{align}
	
	From the above,  using $|f|_2^2\leq |f|_4^2V(\Omega)^{\frac12}$,  we obtain
	\begin{align}
		(SN_2(\Omega))^{-2}|f|_4^4&\leq (SN_2(\Omega))^{-2}|f|_2^4V(\Omega)^{-1}+ 4|f|_2^2\cdot |Df|_2^2= |f|_2^2\cdot \Big\{(SN_2(\Omega))^{-2}V(\Omega)^{-1}|f|_2^2+ 4 |Df|_2^2\Big\} \nonumber\\
		&\leq V(\Omega)^{\frac12}|f|_4^2\cdot \Big\{(SN_2(\Omega))^{-2}V(\Omega)^{-1}|f|_2^2+ 4 |Df|_2^2\Big\} \nonumber.  
	\end{align}
	Simplifying the above yields the conclusion.
}
\qed

\begin{lemma}\label{prop C0 bound by integral}
	{Assume $-\Delta f= \mu\cdot f$ in $\Omega$,  with $\frac{\partial f}{\partial\vec{n}}\big|_{\partial\Omega}= 0$,  then
		\begin{align}
			|f|_\infty\leq |f|_2\cdot C\cdot (1+ \mu)V(\Omega)^{-\frac12},  \nonumber 
		\end{align}
		where $C> 0$ is some universal constant. 
	}
\end{lemma}

\pf
{Applying Lemma \ref{lem Sobolev improved deri ests} on $f^k$ for any $k\geq 1$,  we get 
	\begin{align}
		V(\Omega)^{\frac12}|f^k|_4^2\leq |f^k|_2^2+ 4(SN_2(\Omega))^2V(\Omega)\cdot |D(f^k)|_2^2.  \nonumber 
	\end{align}
	
	Note 
	\begin{align}
		|D(f^k)|_2^2&= \int_\Omega k^2 f^{2k- 2}Df\cdot Df= \frac{k^2}{2k- 1}\int_\Omega D(f^{2k- 1})\cdot Df= \frac{-k^2}{2k- 1}\int_\Omega f^{2k- 1}\cdot \Delta f \nonumber \\
		&= \frac{\mu\cdot k^2}{2k- 1}|f|_{2k}^{2k}.  \nonumber
	\end{align}
	
	The above implies
	\begin{align}
		V(\Omega)^{\frac12}|f|_{4k}^{2k}\leq |f|_{2k}^{2k}\Big\{1+ 4(SN_2(\Omega))^2\cdot \frac{\mu k^2}{2k- 1}V(\Omega)\Big\}.  \nonumber 
	\end{align}
	
	Let $\alpha= 4(SN_2(\Omega))^2$,  therefore we obtain
	\begin{align}
		V(\Omega)^{-\frac{1}{4k}}|f|_{4k}\leq V(\Omega)^{-\frac{1}{2k}}|f|_{2k}\cdot (1+ \alpha\cdot \frac{k^2\mu}{2k- 1}V(\Omega))^{\frac{1}{2k}}.  \nonumber 
	\end{align}
	
	Choose $k= 2^{i- 1}$ for $i= 1, 2, \cdots$,  by induction,  using $\mathcal{D}(\Omega)= 1$ and Theorem \ref{SN2 upper bound},  we have
	\begin{align*}
		|f|_\infty =& \lim_{i\rightarrow\infty}V(\Omega)^{-\frac{1}{2^{i}}}|f|_{2^{i}}\leq V(\Omega)^{-\frac12}|f|_2\cdot \sqrt{\prod_{j= 0}^\infty [1+ \alpha\cdot \frac{2^{2j}\mu}{2^{j+ 1}- 1}]^{2^{-j}}} \\
		\leq& V(\Omega)^{-\frac12}\cdot |f|_2\cdot (1+ \mu)2\alpha\\
		\leq&C(1+\mu)V(\Omega)^{-\frac12}|f|_2.  
	\end{align*}
	
}
\qed

\section{Comparison of first eigenvalue in $2$-dm case}

In this section, we assume $\Omega$ is a convex domain in $\mathbb{R}^2$, with diameter $\mathcal{D}(\Omega)=1$, and width $\mathcal{W}(\Omega)=\epsilon$ for some $\epsilon<\frac{1}{40}$. By Lemma \ref{lem put convex domain into Rn}, we further assume $\Omega\subseteq(0,1)\times(-\epsilon,\epsilon)\subseteq \mathbb{R}^2$, and $(0,0),(1,0)\in\overline{\Omega}$.

 Assume  $\Omega_t=\Omega\cap\{x=t\}=\{(t,y)\in\mathbb{R}^2|y\in(h_-(t),h_+(t))\}$.   Let $h(t)=h_+(t)-h_-(t)$.

Let $u$ be the first Neumann eigenfunction of $\Omega$ with respect to $\mu$, $-\Delta u=\mu\cdot u$, assume $\sup|u|=1$.
For $x\in[0,1]$, define $\bar{u}(x)=\frac{1}{h(x)}\int_{\Omega_x}u(x,y)dy$. We define the modified cross-sectional average $\tilde{u}$ and the error term $\eta$ as follows:
\begin{align}
\hat{u}=\left\{\begin{array}{lll}
&\bar{u}(\epsilon),\quad &x\in [0,\epsilon]; \\
&\bar{u}(x),\quad &x\in[\epsilon,1-\epsilon]; \\
&\bar{u}(1-\epsilon),\quad &x\in[1-\epsilon,1].
\end{array}\right.
\end{align}

\begin{align}\label{tilde u}
c_1=\frac{\int_0^1h\hat{u}}{\int_0^1h},\quad \tilde{u}=\hat{u}-c_1.
\end{align}

\begin{align}
\eta(x)=h\tilde{u}'(x)+\mu\int_0^xh\tilde{u}.
\end{align}

From $\sup|u|=1$, we get $\sup|\bar{u}|\leq1$ and $\sup|\hat{u}|\leq\sup|\bar{u}|\leq1$. So by \eqref{tilde u}, we get $\sup|\tilde{u}|\leq2$.

From \cite[Proposition 2]{Kr}, we get $\mu_1\leq 4j_{0,1}^2\leq C$.

By Li-Yau's gradient estimate \cite{LY}, we have
\begin{align}
|Du|\leq \sqrt{1-u^2}\sqrt{\mu_1}\leq C.
\end{align}

For any $s\in[0,\epsilon]$, we have 
\begin{align}\label{end control}
|\bar{u}(s)-\bar{u}(\epsilon)|\leq|\min_{(x,y)\in\Omega,x\leq\epsilon}u(x,y)-\max_{(x,y)\in\Omega,x\leq\epsilon}u(x,y)|\leq 3\epsilon\cdot\sup|Du|\leq C_1\epsilon.
\end{align}
Same argument for $s\in[1-\epsilon,1]$, we get
\begin{align}
|\bar{u}(s)-\bar{u}(1-\epsilon)|\leq C_1\epsilon.
\end{align}
Assume $h(t_0)=\sup h$, since $\mathcal{W}(\Omega)=\epsilon$, we get $h(t_0)\geq\epsilon$. Since $\Omega$ is convex, we have $\Omega$ contain the quadrilateral formed by $(0,0),(t_0,h_-(t_0)),(1,0),(t_0,h_+(t_0))$. So we get
\begin{align}\label{vol control}
\int_0^1h=V(\Omega)\geq \frac12(h_+(t_0)-h_-(t_0))= \frac12h(t_0)\geq\frac12\epsilon.
\end{align}
So we have
\begin{align}\label{c1 control}
c_1&=\frac{\int_0^1h\hat{u}}{\int_0^1h}=\frac{\int_0^1h(\hat{u}-\bar{u})}{\int_0^1h}\nonumber\\
&=\frac{\int_0^\epsilon h(s)(\bar{u}(\epsilon)-\bar{u}(s))ds+\int_{1-\epsilon}^1h(s)(\bar{u}(\epsilon)-\bar{u}(s))ds}{\int_0^1h}\nonumber\\
&\leq\frac{4C_1\epsilon^3}{\frac12\epsilon}\nonumber\\
&\leq C_2\epsilon^2.
\end{align}

\begin{lemma}\label{h control}
For $x\in(0,1)$, we have
\begin{align}
|h_+'(x)|\leq \frac{h_+(x)}{\min\{x,1-x\}},\quad   \quad|h_-'(x)|\leq \frac{|h_-(x)|}{\min\{x,1-x\}}.
\end{align}
In particular, for $x\in[\epsilon,1-\epsilon]$, we have $|h_+'|,|h_-'|\leq1$.
\end{lemma}

\pf
Since $\Omega$ is convex, we get $h_+$ is a concave function on $[0,1]$, so $h_+'(x)$ is a decreasing function on $[0,1]$.
For any $x\in(0,1)$, by mean value theorem, there exists $x_1\in(0,x)$, such that 
\begin{align}
h_+'(x_1)=\frac{h_+(x)-h_+(0)}{x}=\frac{h_+(x)}{x}, \nonumber
\end{align}
so we get 
\begin{align}
h_+'(x)\leq h_+'(x_1)=\frac{h_+(x)}{x}. \nonumber
\end{align}

Similarly, we get
\begin{align}
h_+'(x)\geq\frac{-h_+(x)}{1-x}. \nonumber
\end{align}

So, we get
\begin{align}
|h_+'(x)|\leq \frac{h_+(x)}{\min\{x,1-x\}} \nonumber
\end{align}
Same argument holds for $-h_-$, so we get the second inequality.
\qed

\begin{lemma}\label{eigenvalue comparison}
We have the following formula:
\begin{align}\label{eigen compare}
\mu=\frac{\int_0^1h(\tilde{u}')^2}{\int_0^1h\tilde{u}^2}-\frac{\int_0^1\eta\tilde{u}'}{\int_0^1h\tilde{u}^2}.
\end{align}
\end{lemma}

\pf
By definition of  $\eta$, we have
\begin{align}
\int_0^1h(\tilde{u}')^2=\int_0^1\eta\tilde{u}'-\mu\int_0^1\tilde{u}'(x)dx\int_0^xh\tilde{u}(s)ds. \nonumber
\end{align}
Since $\int_0^1h\tilde{u}=0$, use integration by parts, we get
\begin{align}
\int_0^1\tilde{u}'(x)dx\int_0^xh\tilde{u}(s)ds=-\int_0^1h\tilde{u}^2. \nonumber
\end{align}
So we get \eqref{eigen compare}
\qed

\begin{lemma}\label{int hu2 control}
	There is a universal constant $C>0$, such that
	\begin{align}
		\int_{-1}^1h\tilde{u}^2\geq C\epsilon.
	\end{align}
\end{lemma}

\pf
By Lemma \ref{prop C0 bound by integral} and \eqref{vol control} we get
\begin{align}
	\int_\Omega u^2\geq C\cdot V(\Omega)\cdot(1+\mu_1)^{-2}\sup|u|^2\geq C\epsilon.
\end{align}

We have 
\begin{align}
	|u(t,y)-\bar{u}(t)|\leq|\sup_yu(t,y)-\inf_yu(t,y)|\leq\sup|Du|\cdot h(t)\leq C\epsilon.
\end{align}
So we get 
\begin{align}
	\bar{u}(t)^2=\fint_{\Omega_t}u(t,y)^2dy-\fint_{\Omega_t}|u(t,y)-\bar{u}(t)|^2dy\geq\fint_{\Omega_t}u^2dy-C\epsilon^2.
\end{align}

So, we have
\begin{align}\label{bar u2 control}
	\int_0^1h\bar{u}^2\geq\int_\Omega u^2-C\int_0^1h\epsilon^2\geq C\epsilon.
\end{align}

By definition of $\tilde{u}$, we get
\begin{align}
	|\int h\tilde{u}^2-\int h\bar{u}^2|&\leq\int_0^{\epsilon}h(t)|(\bar{u}(\epsilon)-c_1)^2-\bar{u}(t)^2|dt
	+\int_{1-\epsilon}^1h(t)|(\bar{u}(1-\epsilon)-c_1)^2-\bar{u}(t)^2|dt\nonumber\\
	&+|\int_{\epsilon}^{1-\epsilon}h\cdot((\bar{u}-c_1)^2-\bar{u}^2)|.
\end{align}
From \eqref{end control} and \eqref{c1 control}, we get
\begin{align}
	|\int h\tilde{u}^2-\int h\bar{u}^2|\leq C\epsilon^3.
\end{align}
Combine this with \ref{bar u2 control}, we get
\begin{align}
	\int_0^1h\tilde{u}^2\geq C\epsilon.
\end{align}
\qed

\begin{lemma}\label{eta control}
There is a universal constant $C>0$, such that
\begin{align}
|\eta(t)|\leq\left\{\begin{array}{l}
C\epsilon^2,\quad\quad\quad\quad\quad\quad\quad t\in[0,\epsilon]\cup[1-\epsilon,1], \\
C\epsilon\cdot h(t)(|h_-(t)'|+|h_+(t)'|)+C\epsilon^3;\quad t\in[\epsilon,1-\epsilon].
\end{array}\right.
\end{align}
\end{lemma}

\pf
From $\Delta u=-\mu_1u$, we get

\begin{align}\label{int by part}
	\mu_1\int_0^th\bar{u}=\mu_1\int_{\Omega\cap\{x\leq t\}}u=-\int_{\Omega_x}u_x(x,y)dy.
\end{align}

First, we compute $\bar{u}'$, we have
\begin{align}\label{bar u'}
\bar{u}'(t)&=\frac{d}{dx}\Big\{\frac{1}{h(t)}\int_{\Omega_t}u(t,y)\Big\} \nonumber\\
&=-\frac{h'}{h^2}\int_{\Omega_t}u+\frac{1}{h}\int_{\Omega_t}u_x+\frac{[u(t,h_+(t))\cdot h_+'(t)-u(t,h_-(t))\cdot h_-'(t)]}{h(t)} 
\end{align}

For $t\in[0,\epsilon]$, we have $\tilde{u}(t)=\bar{u}(\epsilon)-c_1$, and $\tilde{u}'=0$, so we get
\begin{align}
|\eta(t)|=\mu_1|\int_0^th\tilde{u}|\leq C\epsilon\sup|h|\cdot\sup|\tilde{u}|\leq C\epsilon^2
\end{align}
Similarly, for $t\in[1-\epsilon,1]$, we have
\begin{align}
|\eta(t)|=\mu_1|\int_t^1h\tilde{u}|\leq C\epsilon^2
\end{align}

For $t\in[\epsilon,1-\epsilon]$, we have $\tilde{u}'=\bar{u}'$, by \eqref{int by part}, \eqref{bar u'} and Lemma \ref{vertical derivative control},we get
\begin{align}
|\eta(t)|&\leq\Big|-\frac{h'}{h}\int_{\Omega_t}udy+[u(t,h_+(t))\cdot h_+'(t)-u(t,h_-(t))\cdot h_-'(t)]\Big|\nonumber\\
&+\Big|\mu_1\int_0^\epsilon h(s)\cdot[\bar{u}(\epsilon)-c_1-\bar{u}(s)]ds\Big|+\Big|\mu_1\int_{\epsilon}^{1-\epsilon}h(s)\cdot c_1ds\Big| \nonumber\\
&\leq\Big|\frac{\int_{\Omega_t}[u(t,h_+(t))-u(t,y)]\cdot h_+'(t)-[u(t,h_-(t))-u(t,y)]\cdot h_-'(t)dy}{h(t)}\Big|+C\epsilon^3\nonumber\\
&\leq[|h_+(t)'|+|h_-(t)'|]\cdot \sup|D_yu(t,y)|\cdot h(t)+C\epsilon^3\nonumber\\
&\leq C\epsilon\cdot h(t)(|h_-(t)'|+|h_+(t)'|)+C\epsilon^3
\end{align}

\qed

\begin{lemma}\label{int eta u control}
	There is a universal constant $C>0$, such that 
	\begin{align}
		\int_0^1\eta\tilde{u}'\leq C\epsilon^3.
	\end{align}
\end{lemma}
\pf
Assume $h(t_0)=\sup h(t)$, since $\mathcal{W}(\Omega)=\epsilon$, we get $h(t_0)\geq\epsilon$. Since $\Omega$ is convex, we get $h'\geq0$ on $[0,t_0]$ and $h'\leq0$ on $[t_0,1]$. Since $\int h\tilde{u}=0$, we get
\begin{align}
\Big|\frac{\int_0^t\mu_1h(s)\tilde{u}(s)ds}{h(t)}\Big|=\left\{\begin{array}{l}
\Big|\frac{\int_0^t\mu_1h(s)\tilde{u}(s)ds}{h(t)}\Big|\leq C,\quad t\in[0,t_0]; \nonumber\\
\Big|\frac{\int_{t}^1\mu_1h(s)\tilde{u}(s)ds}{h(t)}\Big|\leq C;\quad t\in[t_0,1].
\end{array}\right.
\end{align}
So
\begin{align}\label{eta tilde u}
|\int_0^1\eta\tilde{u}'|&=|\int_{\epsilon}^{1-\epsilon}\eta\cdot(\frac{\eta}{h}-\frac{1}{h}\int_{-1}^t\mu_1h\tilde{u})dt|\nonumber\\
&\leq\int_{\epsilon}^{1-\epsilon}(\frac{\eta^2}{h}+C|\eta|)\nonumber\\
&\leq \int_{\epsilon}^{1-\epsilon}[C\epsilon^2h+C\frac{\epsilon^6}{h}+C\epsilon h\cdot(|h_+'|+|h_-'|)]+C\epsilon^3.
\end{align}
Since $h$ is convex, we get for $t\in[\epsilon,t_0]$, we have
\begin{align}
h(t)\geq \frac{t}{t_0}h(t_0)+\frac{t_0-t}{t_0}h(0)\geq\epsilon h(t_0)\geq\epsilon^2,
\end{align}
for $t\in[t_0,1-\epsilon]$, we have similar argument, so we get for $t\in[\epsilon,1-\epsilon]$, we have
\begin{align}\label{h bound}
h(t)\geq\epsilon^2.
\end{align}
Since $h(t)=h_+(t)-h_-(t)$, we get 
\begin{align}\label{h' bound}
\int_{\epsilon}^{1-\epsilon}h_+|h_+'|\leq\int_0^1\frac{1}{2}|(h_+^2)'|=\sup|h_+^2|\leq\epsilon^2.
\end{align}
Similar argument holds for $h_-$.
Use \eqref{h bound} and \eqref{h' bound} in \eqref{eta tilde u}, we get
\begin{align}
|\int_0^1\eta\tilde{u}'|\leq C\epsilon^3.
\end{align}
\qed

\begin{theorem}\label{eigen control}
	Let $\Omega\subset\mathbb{R}^2$ be a convex domain with diam $\mathcal{D}(\Omega)=1$, and width $\mathcal{W}(\Omega)=\epsilon$. Then, there is a universal constant $C>0$, such that
	\begin{align}
		\mu_1(N)\geq\mu_1(\Omega)\geq\mu_1(N)-C\epsilon^2.
	\end{align}
\end{theorem}

\pf
The left hand side of the inequality follows from Proposition \ref{prop upper bound of mu}.

For the right hand side of the inequality, we additionally let constant $C>6400 j_{0,1}^2$, if width $\mathcal{W}(\Omega)=\epsilon\geq\frac{1}{40}$, since $\mu_1(N)\leq 4j_{0,1}^2$, we get $\mu_1(N)-C\epsilon^2\leq0\leq\mu_1(\Omega)$.

If $\epsilon<\frac{1}{40}$, from Lemma \ref{eigenvalue comparison}, Lemma \ref{int hu2 control} and Lemma \ref{int eta u control}, we get
\begin{align}
	\mu_1(\Omega)=\mu=&\frac{\int_0^1h(\tilde{u}')^2}{\int_0^1h\tilde{u}^2}-\frac{\int_0^1\eta\tilde{u}'}{\int_0^1h\tilde{u}^2}\\
	\geq&\mu_1(N)-C\epsilon^2.
\end{align}
\qed

\section{The eigenvalue estimate of $N$}

In this section, we give another proof for the simplicity of eigenvalue on convex domain in $\mathbb{R}^2$ by estimate the eigenvalue of collapsing segment $N$.

\begin{lemma}\label{lem ln concave for slice of convex domain}
{Let $\Omega\subset\mathbb{R}^n$ be a convex set with smooth boundary, and $H(x)$ be the $(n-1)$-dim Hausdorff measure of $\Omega_x\vcentcolon=\{y\in\mathbb{R}^{n-1}:(x,y)\in\Omega\}$, then we have $$(H^{\frac{1}{n- 1}})''\leq 0$$.
}
\end{lemma}

\pf
We only need to show $H^\frac{1}{n-1}$ is a concave function.

For any $-1\leq a<b\leq1$, we have $\Omega_a=\{y\in\mathbb{R}^{n-1}:(a,y)\in\Omega\}$,  and $\Omega_b=\{y\in\mathbb{R}^{n-1}:(b,y)\in\Omega\}$ are bounded subsets in $\mathbb{R}^{n-1}$. Since $\Omega$ is convex, we get 
\begin{align}
\frac{\Omega_a+\Omega_b}{2}\subset\Omega_{\frac{a+b}{2}}, \nonumber
\end{align}
where $\frac{\Omega_a+\Omega_b}{2}=\{\frac{p+q}{2}, \,p\in\Omega_a,q\in\Omega_b\}$.

By \cite[Theorem $3.2.41$]{Federer}, we get  
\begin{align}
\mathcal{H}^{n-1}(\Omega_a)^\frac{1}{n-1}+\mathcal{H}^{n-1}(\Omega_b)^\frac{1}{n-1}\leq\mathcal{H}^{n-1}(\Omega_a+\Omega_b)^\frac{1}{n-1}=2\mathcal{H}^{n-1}(\frac{\Omega_a+\Omega_b}{2})^\frac{1}{n-1} \nonumber
\end{align}
so we get 
\begin{align}
\frac{\{\mathcal{H}^{n-1}(\Omega_a)^\frac{1}{n-1}+\mathcal{H}^{n-1}(\Omega_b)^\frac{1}{n-1}\}}{2}\leq\mathcal{H}^{n-1}(\Omega_\frac{a+b}{2})^\frac{1}{n-1} \nonumber
\end{align}

Since 
\begin{align}
H(x)=\mathcal{H}^{n-1}(\Omega_x)\times(\int_{-1}^1\mathcal{H}^{n-1}(\Omega_t)dt)^{-1},\quad x\in[-1,1]
\end{align}
we get 
\begin{align}
\frac{\Big(H(a)^\frac{1}{n-1}+H(b)^\frac{1}{n-1}\Big)}{2}\leq \frac{H(\frac{a+b}{2})^\frac{1}{n-1}}{2} \nonumber
\end{align}
So $H^\frac{1}{n-1}$ is concave function.
\qed

The eigenvalue of $N$ is $\lambda$ with respect to
\begin{align}
\psi''(t)+ (\ln H)' \psi'(t)+ \lambda\cdot \psi(t)= 0, \quad \quad \quad \forall t\in [0, 1], \psi'(0)=\psi'(1)= 0. \nonumber 
\end{align}

Define $V= \frac{3}{4}[(\ln H)']^2- \frac{H''}{2H}$ and $L_2(w)= w''- V\cdot w$, where $w\in C_0^\infty[0, 1]$. 

\begin{lemma}\label{lem eigen of N and with potential}
{For any $k\in \mathbb{Z}^+$, we have 
\begin{align}
\lambda_k(L_2)= \mu_{k}(N). \nonumber 
\end{align}
}
\end{lemma}

\pf
{If $v$ minimizes the quotient $\frac{\int_0^1 h(t)(\phi'(t))^2dt}{\int_0^1 h(t)\phi(t)^2dt}$ among functions $\phi$ satisfying the assumption, then for any function $f$ satisfying $\int_0^1 fhdt= 0$, we have 
\begin{align}
0= \frac{d}{ds}\big|_{s= 0}\Big(\frac{\int_0^1 h(t)((v+ sf)'(t))^2dt}{\int_0^1 h(t)(v+ sf)(t)^2dt}\Big). \nonumber 
\end{align}
Simplifying the above equation, we get 
\begin{align}
\int_0^1 hv'f'- \lambda hvfdt= 0, \label{first variation vanishes for restrcted f}
\end{align}
where $\lambda$ is the minimum value of the quotient.

Note for $f\equiv 1$, (\ref{first variation vanishes for restrcted f}) also holds because of $\int_0^1 hvdt= 0$.

Note $1= \alpha\cdot h+ h^\perp$, where $\int_0^1 h\cdot h^\perp= 0, \alpha\in \mathbb{R}$. If $\alpha= 0$, then $\int_0^1 h=0$, however $h\geq 0$, hence $h\equiv 0$, which is the contradiction. Now from $\alpha\neq 0$, we know (\ref{first variation vanishes for restrcted f}) also holds for $f= h$. Hence (\ref{first variation vanishes for restrcted f}) holds for any $f$.

From (\ref{first variation vanishes for restrcted f}), integration by parts yields
\begin{align}
(hv'f)\big|_0^1- \int_0^1 \Big\{(hv')'+ \lambda hv\Big\}\cdot fdt= 0. \label{with boundary integral}
\end{align}
Choose $f\in C_0^\infty[0, 1]$ in the above equation, we get $\int_0^1 \Big\{(hv')'+ \lambda hv\Big\}\cdot fdt= 0$ for any such $f$. Therefore, we get 
\begin{align}
(hv')'+ \lambda hv= 0. \label{2nd order ODE equa}
\end{align}

Plugging (\ref{2nd order ODE equa}) into (\ref{with boundary integral}), we obtain $(hv'f)\big|_0^1= 0$ for any $f$, hence we have $hv'(0)= hv'(1)= 0$. 

Let $w= \sqrt{h} v'$, then we have
\begin{align}
w''= -\lambda w+ (\frac{3}{4}(\frac{h'}{h})^2- \frac{h''}{2h})w, \quad \quad and \quad \quad w(0)= w(1)= 0. \nonumber
\end{align}
}
\qed

Then we have
\begin{prop}\label{prop ODE eigenv bound by without potential}
{For any $k\in \mathbb{Z}^+$, we have $\displaystyle \mu_{k}(N)\geq k^2\pi^2$. 
}
\end{prop}

\pf
{\textbf{Step (1)}. Define $L_1(w)= w''$, where $w\in C_0^\infty[0, 1]$, then $\lambda_k(L_1)= k^2\pi^2$. 

Note $V= (\frac{3}{4}(\frac{h'}{h})^2- \frac{h''}{2h})\geq -\frac{1}{2}(\ln h)''$, then from the assumption we get $V\geq 0$. Then we get 
\begin{align}
\frac{\int_0^1 w(-L_2(w))dt}{\int_0^1 w^2dt} \geq \frac{\int_0^1 w(t)\cdot (-L_1(w))dt}{\int_0^1 w(t)^2dt}, \quad \quad \forall w\in  C_0^\infty[0, 1].\label{Rayleigh quot compar}
\end{align}

Now from (\ref{Rayleigh quot compar}), we obtain
\begin{align}
\lambda_k(L_1)= \inf_{V_k\subseteq C_0^\infty[0, 1]}\sup_{w\in V_k- \{0\}} \frac{\int_0^1 w(t)\cdot (-L_1(w))dt}{\int_0^1 w(t)^2dt}\leq \inf_{V_k\subseteq C_0^\infty[0, 1]}\sup_{w\in V_k- \{0\}} \frac{\int_0^1 w(-L_2(w))dt}{\int_0^1 w^2dt} = \lambda_k(L_2)=\mu_k(N), \nonumber 
\end{align}
where $V_k$ is any $k$-dim subspace. 
}
\qed

Based on the above, we get another proof for the simplicity of eigenvalue on convex domain in $\mathbb{R}^2$,
\begin{cor}
Let $\Omega$ be a convex domain in $\mathbb{R}^2$ with diameter 1 and width $\epsilon$. Where 
$$
\epsilon<\frac{2(\pi-j_{0,1})(k\pi+j_{0,1})}{(k+1)^2\pi^2[(k+1)^2\pi^2+1]}
$$
then the first $k$ Neumann eigenvalue of $\Omega$ is simple.
\end{cor}

\pf
Let $N$ be the corresponding collapsing segment of $\Omega$, from Proposition \ref{prop ODE eigenv bound by without potential} and \cite{Kr} we have for any integer $k>0$,
\begin{align}
k^2\pi^2\leq \mu_k(N)\leq (2j_{0,1}+(k-1)\pi)^2, \quad \quad \quad \forall k\in \mathbb{Z}^+. \nonumber
\end{align}

Applying Theorem \ref{thm Omegai and N comparsion}, we get $\mu_k(\Omega)\leq\mu_k(N)< (1-2\epsilon(1+\mu_{k+1}(N)))\mu_{k+1}(N)\leq\mu_{k+1}(\Omega)$, for $s=1,2,\cdots,k$. So the first $k$ eigenvalue of $\Omega$ is simple.
\qed

\begin{remark}
In the above proof, the bound for width $\epsilon$ is $\epsilon\sim k^{-3}$, which is weaker than Theorem \ref{2dim eigen simple} with $\epsilon\sim k^{-1}$. The main problem is that we do not have a direct bound for the gap between eigenvalue of $N$. In the above proof, we use the fact that the upper bound of eigenvalue of $N$ is smaller than the lower bound of the next eigenvalue of $N$. This no longer holds for higher dimension. 
\end{remark}

\appendix
\section{The estimate of Dirichlet-Neumann heat kernel}
\begin{center}
    By Guoyi Xu\footnote{Department of Mathematical Sciences, Tsinghua University, Beijing, P. R. China. E-mail address: guoyixu@tsinghua.edu.cn}
\end{center}

The method of using Dirichlet-Neumann heat kernel to estimate the Neumann eigenvalue in this section,  is inspired by the argument of the lower bound for heat kernel in \cite{CY} and \cite{LY-parabolic}.  

\begin{lemma}\label{lem solution to Robin para equ}
	{Assume $f\in C^\infty[0, \rho_1]$ with 
		\begin{align}
			f'(0)= f(\rho_1)=0, \quad \quad \quad f\geq 0, \quad \quad \quad f'(y)\leq 0, \quad \quad \forall y\in [0, \rho_1]. \label{special f}
		\end{align}
		Furthermore, assume $v(y, t): [0, \rho_1]\times [0, \infty)\rightarrow \mathbb{R}$ is a solution to
		\begin{equation}\label{parabolic equa}
			\left\{
			\begin{array}{rl}
				&v_t- v_{yy}= 0\ ,  \quad \quad \quad \quad\quad  \quad\quad \quad \forall (y, t)\in (0, \rho_1)\times (0, \infty)\\
				&v(y, 0)= f(y), \quad \quad  \forall y\in [0, \rho_1], \quad \text{and} \quad v_y(0, t)= v(\rho_1, t)= 0 \ ,  \quad \quad \quad \quad \forall t\in [0, \infty)
			\end{array} \right.
		\end{equation}
		Then $v_y(y, t)\leq 0$ for all $t\geq 0, y\in [0, \rho_1]$.
	}
\end{lemma}

\pf
{From Hopf's strong Maximum principle and the Maximum principle,  we get
	\begin{align}
		\min_{(y, t)\in [0, \rho_1]\times [0, \infty)}v(y, t)\geq \min_{(x, t)\in ([0, \rho_1]\times \{0\}) \cup(\{0, \rho_1\}\times \mathbb{R}^+)}v(x, t)= \min_{(x, t)\in ([0, \rho_1]\times \{0\}) \cup(\{\rho_1\}\times \mathbb{R}^+) }v(x, t)\geq 0.  
	\end{align}
	This implies $\frac{\partial}{\partial y}v(\rho_1, t)\leq 0$ because $v(\rho_1, t)= 0$.  
	
	Let $u= v_y$, then $u_t- u_{yy}= 0$.   From the assumption and the above,  we have
	\begin{align}
		u(y,  0)= f'(y)\leq 0,  \quad \quad u(0, t)= v_y(0,  t)= 0,  \quad \quad u(\rho_1, t)= v_y(\rho_1, t)\leq 0.  \nonumber 
	\end{align}
	Applying the Maximum principle on $u$,  we get $u(y, t)\leq 0$,  and the conclusion follows.
}
\qed

\begin{definition}\label{def DN eigenfunction and heat kernel}
	{For a bounded domain $\Omega\subseteq \mathbb{R}^n$ (where $n\geq 2$) with $\partial \Omega= \gamma_1\cup \gamma_2$, where $\gamma_1$ is Lipschitz and $\gamma_2$ is a convex boundary (or $C^2$ boundary); the \textbf{Dirichlet-Neumann heat kernel} $\mathcal{K}(x, y, t)$ with respect to $(\Omega, \gamma_1, \gamma_2)$,  is the function satisfying the follows:
		\begin{align}
			&\mathcal{K}(x,y, t)= \mathcal{K}(y, x, t),  \quad \quad \mathcal{K}(x,\cdot ,  t)\big|_{\gamma_1}= \frac{\partial}{\partial \vec{n}}\mathcal{K}(x,\cdot ,  t)\big|_{\gamma_2}= 0,  \nonumber \\
			& (\frac{\partial}{\partial t}- \Delta_y)\mathcal{K}(x,y, t)= 0,  \quad \quad \lim_{t\rightarrow 0}\int_\Omega \mathcal{K}(x, y, t)f(x)dx= f(y),  \quad \quad \forall f\in L^\infty(\Omega).  \nonumber 
		\end{align}
	}
\end{definition}

\begin{definition}\label{def DN eigenfunc}
	{We say $\varphi$ is a \textbf{Dirichlet-Neumann eigenfunction with respect to $\lambda$ and $(\Omega, \gamma_1, \gamma_2)$}, if
		\begin{align}
			\Delta \varphi= -\lambda\cdot \varphi,  \quad \quad \varphi\big|_{\gamma_1}=0,  \quad \quad \frac{\partial}{\partial\vec{n}}\varphi\big|_{\gamma_2}= 0.  \nonumber 
		\end{align}
	}
\end{definition}

For any Dirichlet-Neumann heat kernel $\mathcal{K}(x, y, t)$ with respect to $(\Omega, \gamma_1, \gamma_2)$, it is known that
\begin{align}
	\mathcal{K}(x, y, t)= \sum_{i= 1}^\infty e^{-\lambda_i\cdot t}\varphi_i(x)\varphi_i(y),  \nonumber 
\end{align}
where $\{\varphi_i\}_{i= 1}^\infty$ is a basis of $L^2(\Omega)$ and $\varphi_i$ is the Dirichlet-Neumann eigenfunction with respect to $\lambda_i$. 

\begin{theorem}\label{thm DN-heat kernel comparison}
	{Assume $\partial\Omega= \Lambda\cup\gamma$ with $\vec{n}\big|_{\Lambda}\cdot \partial y_2\leq 0$, where $\vec{n}$ is the outward unit normal to $\partial{\Omega}$; and $(0, 0)\in \Omega\subseteq \{(y_1, y_2)\in \mathbb{R}^n: y_1\in \mathbb{R}^{n-1}, 0\leq y_2\leq \rho_1\}$. Then the first Dirichlet-Neumann eigenvalue $\lambda_1(\Omega)$ with respect to $(\Omega, \gamma, \Lambda)$ satisfies the following inequality:
		\begin{align}
			\lambda_1(\Omega)\geq \frac{\pi^2}{4\rho_1^2}. \nonumber 
		\end{align}
	}
\end{theorem}

\pf
{\textbf{Step (1)}.  Let $\mathcal{K}(x, y, t)$ be the Dirichlet-Neumann heat kernel with respect to $(\Omega,  \gamma,  \Lambda)$, and $\mathcal{H}(x, y, t)$ be the Dirichlet-Neumann heat kernel with respect to $([0, \rho_1],  \rho_1,  0)$.
	
	For any $f$ satisfying (\ref{special f}), we define 
	\begin{align}
		F(y, t)\vcentcolon= \int_{\Omega} \mathcal{K}((z_1, z_2), y, t)f(z_2)dz_1dz_2- \int_0^{\rho_1} \mathcal{H}(z, y_2, t)f(z)dz, \nonumber 
	\end{align}
	where $y= (y_1, y_2)\in \Omega$. 
	
	Direct computation yields
	\begin{align}
		(\frac{\partial}{\partial t}- \Delta_y)F(y, t)= 0, \quad \quad \quad \forall (y, t)\in \Omega\times (0, \infty). \nonumber 
	\end{align}
	
	From the Maximum Principle for heat equation, we get 
	\begin{align}
		\max_{(x, t)\in \Omega\times [0, \infty)} F(x, t)\leq \max_{(y, t)\in (\partial \Omega\times \mathbb{R}^+)\cup (\Omega\times \{0\})} F(y, t). \label{direct by max prin}
	\end{align}
	
	\textbf{Step (2)}. For $y\in \Lambda$, from the assumption $\vec{n}\cdot \partial y_2\leq 0$ and Lemma \ref{lem solution to Robin para equ}, we have 
	\begin{align}
		\frac{\partial F}{\partial\vec{n}}(y, t)= 0- (\vec{n}\cdot \partial y_2)\cdot  \frac{d}{dy_2}\int_0^{\rho_1} \mathcal{H}(z, y_2, t)f(z)dz\leq 0. \nonumber 
	\end{align}
	Therefore from the Hopf's Strong Maximum Principle and (\ref{direct by max prin}), we get that 
	\begin{align}
		\max_{(x, t)\in \Omega\times [0, \infty)} F(x, t)\leq \max_{(y, t)\in (\gamma \times \mathbb{R}^+)\cup (\Omega\times \{0\})} F(y, t). \nonumber 
	\end{align}
	
	Note $F\big|_{ (\gamma \times \mathbb{R}^+)\cup (\Omega\times \{0\})}\leq 0$ by the Dirichlet-Neumann boundary condition of $\mathcal{K}$ and $\mathcal{H}$, we get 
	\begin{align}
		\max_{(x, t)\in \Omega\times [0, \infty)} F(x, t)\leq 0. \nonumber 
	\end{align}
	
	\textbf{Step (3)}. By the above, we know
	\begin{align}
		\int_{\Omega} \mathcal{K}((z_1, z_2), y, t)f(z_2)dz_1dz_2\leq \int_0^{\rho_1} \mathcal{H}(z, y_2, t)f(z)dz. \nonumber 
	\end{align}
	
	Using the representation formula of Dirichlet-Neumann heat kernel by eigenfunctions,  we get 
	\begin{align}
		\sum_{i= 1}^\infty \int_{\Omega} e^{-\lambda_i t}\varphi_i(z_1, z_2)\cdot \varphi_i(y) f(z_2)dz_1dz_2\leq \sum_{i= 1}^\infty \int_0^{\rho_1} e^{-\tilde{\lambda_i t}}\tilde{\varphi}_i(z)\tilde{\varphi}_i(y_2)f(z)dz, \nonumber 
	\end{align}
	where $\lambda_i$ is the $i$-th Dirichlet-Neumann eigenvalue with respect to $(\Omega, \gamma, \Lambda)$,  and similar for $\tilde{\lambda}_i$ with respect to $([0, \rho_1], \rho_1, 0)$.  
	
	The $1$st Dirichlet-Neumann eigenvalue $\tilde{\lambda}_1$ with respect to $([0, \rho_1], \rho_1, 0)$ is equal to $\frac{\pi^2}{4\rho_1^2}$.
	
	Let $t\rightarrow \infty$,  note $\varphi_1> 0$ and $\tilde{\varphi}_1> 0$,  we get $\lambda_1\geq \tilde{\lambda_1}= \frac{\pi^2}{4\rho_1^2}$.      
}
\qed

\end{document}